\documentclass[11pt, english]{article}
\usepackage[utf8]{inputenc}
\usepackage[T1]{fontenc}

\usepackage{graphicx,epsfig,color}
\usepackage{booktabs}
\usepackage{amsmath, amssymb, amsthm, amsfonts}
\usepackage{mathrsfs}
\usepackage{mathtools}
\usepackage{bm}
\usepackage{epstopdf}
\usepackage{ifthen}
\usepackage{bbm}
\usepackage{xfrac}
\graphicspath{{eps/}}
\usepackage{bookmark}
\usepackage{hyperref}
\usepackage{comment}
\usepackage{marginnote}
\usepackage{clchanges}

\newcommand{\R}{ \mathbb{R}}
\newcommand{\Z}{ \mathbb{Z}}
\newcommand{\T}{\mathbb{T}}

\newcommand{\dint}{{\rm d}}

\DeclareMathOperator*{\argmin}{arg\,min}

\DeclareMathOperator{\Conv}{conv}
\DeclareMathOperator{\Cone}{cone}

\DeclareMathOperator{\intr}{int}
\DeclareMathOperator{\diam}{diam}
\DeclareMathOperator{\CPWL}{CPWL}
\DeclareMathOperator{\CH}{CH}
\DeclareMathOperator{\adj}{adj}
\DeclareMathOperator{\St}{St}
\DeclareMathOperator{\Lip}{Lip}
\DeclareMathOperator{\HLip}{HLip}
\DeclareMathOperator{\HTV}{HTV}
\DeclareMathOperator{\TV}{TV}

\DeclareMathOperator{\spn}{span}
\DeclareMathOperator{\dist}{dist}

\def\V#1{{\boldsymbol{#1}}}
\def\Spc#1{{\mathcal{#1}}}
\def\M#1{{\bf{#1}}}
\def\Op#1{{\mathrm{#1}}}

\newtheorem{theorem}{Theorem}[section]
\newtheorem{definition}[theorem]{Definition}
\newtheorem{proposition}[theorem]{Proposition}

\newtheorem{lemma}[theorem]{Lemma}

\newtheorem{example}[theorem]{Example}

\newcommand{\Tr}{\mathsf{T}}

\title{
Quasi-norms on $\R^d$: Characterization and Approximation by Conic Splines
\thanks{The research leading to these results has received funding from the Swiss National Science Foundation (SNSF) under Grant 200020\_219356.
}}

\author{
Michael Unser\thanks{Biomedical Imaging Group, \'Ecole polytechnique f\'ed\'erale de Lausanne (EPFL),
Station 17, CH-1015, Lausanne, Switzerland ({\tt michael.unser@epfl.ch}). }
 }

\begin{document}

\maketitle

\begin{abstract}
We first show that a regularizer on $\mathbb{R}^d$ that is continuous, symmetric, and amplitude-equivariant is of the form $g=q^\gamma$ with $\gamma>0$, where $q$ is continuous and absolutely homogeneous; we call the set of such $q$ the CH class, whose convex members are the seminorms. Within the CH class, we single out two configurations of interest for regularization: the quasi-norms, which are exactly the positive-definite members of the family, and the quasi-seminorms with null space control, which are the composition of a quasi-norm with a linear map that sets the null space. 
We then propose a universal parametrization of the class by conic splines, the continuous functions that are linear on each cone of a polyhedral fan. Specifically, we show that the polyhedral CH functionals are exactly the symmetric non-negative conic splines, and that these admit a stable nodal representation in a basis of conic B-splines. Finally, we determine the rates of approximation: the conic spline interpolant converges at the second-order rate $h^2\asymp n^{-2/(d-1)}$, in the $L_\infty$-norm on the hyper-sphere when the angular profile is twice differentiable, and in the $L_1$-norm when its spherical Hessian variation is bounded, which covers the polyhedral members of the class. The rate is sharp, already for the Euclidean norm. The planar case is worked out in full, with the trigonometric hat functions as basis and an $\ell_1$-norm of second differences as regularity index.
\end{abstract}

\pagebreak
\tableofcontents
\pagebreak

\section{Introduction}
\label{Sec:Introduction}

The functionals that this paper studies arise from the {\em learning} of a regularizer. The reconstruction of a signal $\V s\in\R^d$ from the noisy linear measurements $\V y=\M H\V s+\V \epsilon\in\R^M$ is formulated as
\begin{align}
\label{Eq:VarReconstruction}
\V s^\star=\argmin_{\V s \in \R^d} \tfrac12\|\V y-\M H\V s\|_2^2+\lambda\, g_{\V \theta}(\V s),
\end{align}
where the regularizer $g_{\V \theta}$ carries the prior information and where its parameters $\V \theta$ are trained on representative data \cite{Benning2018,McCann2019,Arridge2019,Mukherjee2023}. In practice, the quality of the reconstruction is tuned through a proper adjustment of the regularization strength $\lambda>0$. The appeal of \eqref{Eq:VarReconstruction}---as opposed to a direct regression by a neural network---is that the reconstruction remains the solution of an optimization problem, so that the existence of a solution, the convergence of the iterative scheme, and the stability of the map $\V y\mapsto\V s^\star$ remain within reach; but those statements are only as strong as the structural constraints that one is willing to impose on $g_{\V \theta}$. Three requirements are natural: {\em amplitude equivariance}, which asks that the amplification of the data by a factor $\alpha\neq0$ amplify the reconstruction by the same factor up to a readjustment of $\lambda$;
{\em symmetry}, $g_{\V \theta}(-\V s)=g_{\V \theta}(\V s)$, since the sign of a contrast carries no prior preference; and {\em continuity}, the bare minimum for any statement about convergence. Under these three requirements {\em together with convexity}, the admissible regularizers collapse to the powers $g=p^\gamma$, $\gamma\ge1$, of a seminorm $p$ \cite[Theorem 2.3]{UnserDucotterd2026}. The same paper also provides a universal parametrization of such functionals by means of learnable polyhedral norms and a demonstration of usage for biomedical image reconstruction.

\subsection{Objectives and Contributions}
The purpose of the present study is to extend the above framework by dropping the convexity constraint, which the sparsity-promoting penalties of practical interest---foremost among them the $\ell_p$ ``norms'' with $0<p<1$---do not satisfy. This raises three questions, which organize the paper: the identification of the class of relevant regularization functionals; the formulation of a parametric framework that is universal for the class; and the rates of approximation that the parametrization achieves.

\medskip
\noindent{\bf\em 1. Identification and categorization of the class.}
We first show that the continuous, symmetric, and amplitude-equivariant functionals are exactly the powers $g=q^\gamma$, $\gamma>0$, of the functionals $q$ that are continuous and absolutely homogeneous,
\begin{align}
\label{Eq:AmplitudeEquivariance}
\forall \V x \in \R^d, \ \forall \alpha \in \R: \quad q(\alpha \V x)=|\alpha| \, q(\V x),
\end{align}
which we refer to as the CH class; in particular, $g$ is convex if and only if $\gamma\ge1$ and $q$ is a seminorm (Proposition~\ref{Prop:AERegularizers}). Since the exponent $\gamma$ is immaterial for \eqref{Eq:VarReconstruction}, the object to be learned is $q$ itself. Within the CH class on $\R^d$, the quasi-norms are exactly the positive-definite members (Theorem~\ref{Theo:QuasiNormAuto}). Likewise, the quasi-seminorms are exactly the members that are comparable to a seminorm; their null space $\Spc N_q$ is then a linear subspace, and $q$ is comparable to the distance to it (Theorem~\ref{Theo:SemiQuasiNorm}).
The null space has an essential role in regularization theory: it is the mechanism by which one favors a type of solution---the constant signals, say, when $\Spc N_q$ is the kernel of a finite-difference operator---while the data alone decide the component of the reconstruction along $\Spc N_q$.
In our framework, this type of regularizer is implemented as the composition $\bar q\circ\M P$ of a quasi-norm $\bar q$ with a matrix $\M P$ that sets the null space.
Theorem~\ref{Theo:QuasiSeminormStructure} characterizes these regularizers as the quasi-seminorms that are invariant along their null space. On the geometric side, the continuous quasi-norms are the gauges of the symmetric star-shaped bodies with a continuous positive radial function (Theorem~\ref{Theo:GaugeSemiQuasiNorm}).

\medskip
\noindent{\bf\em 2. Universal parametrization by conic splines.}
A parametrization of the class should be: (i) universal, in that it covers the class in the limit of infinitely many parameters; (ii) stable, so that the parameters can be adjusted by a descent algorithm; and (iii) free of constraints, so that the admissibility of the learned functional is guaranteed by construction, as the polyhedral norms achieve in the convex theory.

To this end, we introduce the {\em conic splines}: the continuous functions that are linear on each cone of a complete simplicial fan of $\R^d$ or, equivalently, the $1$-homogeneous extensions of the piecewise-linear functions on a triangulation of the hyper-sphere. We show that the symmetric non-negative conic splines are {\em precisely} the polyhedral members of the CH class (Theorem~\ref{Theo:PolyQuasiSemiNorms}), and that the {\em conic B-splines}---the conic analogues of the Courant hat functions---form a nodal basis of the conic splines on a fan, with the values of the functional along the rays as coefficients and with explicit Riesz bounds (Proposition~\ref{Prop:ConicBsplines}). The parametrization of a polyhedral CH functional by its nodal values is thus free of constraints up to the non-negativity and the symmetry of those values. By Theorem~\ref{Theo:QuasiNormAuto}, the parametrized functional is a polyhedral quasi-norm if and only if all the nodal values are (strictly) positive. The composition of the latter with a matrix $\M P\in\R^{d'\times d}$ then yields a constraint-free parametrization of the quasi-seminorms with the prescribed null space $\Spc N_q=\ker(\M P)$.

\medskip
\noindent{\bf\em 3. Rates of approximation and their sharpness.}
Universality is a qualitative statement, which, for approximation purposes, needs to be complemented by a rate:
how many rays a prescribed accuracy costs, which regularity of the profile this depends on, and whether the rate is the best possible.

For a $1$-homogeneous $q$ that is $C^2$ away from the origin, we show that the conic spline interpolant $q_n=\sum_i q(\M r_i)\varphi_i$ on a shape-regular fan with rays $\M r_1,\ldots,\M r_n$ and maximal spherical cell diameter $h_n$ satisfies
\begin{align*}
\|q-q_n\|_{L_\infty(\mathbb S^{d-1})}\le C\,h_n^2\,\HLip_{\mathbb S^{d-1}}(q),
\end{align*}
where $\HLip_{\mathbb S^{d-1}}(q)$ is the supremum of the spectral norm of the Hessian of $q$ over the hyper-sphere; for quasi-uniform rays, $h_n\asymp n^{-1/(d-1)}$ and the rate is $n^{-2/(d-1)}$ (Theorem~\ref{Theo:ConicSplineApprox}). The rate is sharp: by Gruber's lower bound on the approximation of a convex body by an inscribed polytope with $n$ vertices \cite{Gruber1993}, no conic spline on a fan with at most $n$ rays approximates the Euclidean norm better than $c_d\,n^{-2/(d-1)}$, so that the approximation saturates at the second order however smooth the quasi-norm is (Theorem~\ref{Theo:ConicSplineApprox}, Item 3). We then extend our second-order estimate to the CH functionals with a bounded {\em spherical Hessian variation} $\HTV_{\mathbb S^{d-1}}(q)$, a quantity defined by duality that is finite for the polyhedral members of the class: a preaveraged conic spline then converges in $L_1(\mathbb S^{d-1})$ at the rate $h_n^2\,\HTV_{\mathbb S^{d-1}}(q)$ (Theorem~\ref{Theo:HTVApprox}). In the planar case, the conic B-splines on a uniform polar grid are trigonometric hat functions and the spherical Hessian variation of a polyhedral profile is an $\ell_1$-norm of a second difference of the nodal values (Section~\ref{Sec:Planar}), which is the penalty of the learning problem.

\subsection{Related Works}
Non-convex sparsity-promoting penalties have a long history in signal reconstruction: the $\ell_p$ quasi-norms with $0<p<1$ recover sparse signals under weaker conditions than the $\ell_1$-norm \cite{Chartrand2007,Foucart2009,FoucartRauhut2013}, the reweighted $\ell_1$ schemes approximate them \cite{Candes2008}, and the non-convex penalties of statistics and image processing \cite{Fan2001,Zhang2010,Nikolova2005,Nikolova2010,Soubies2015,Selesnick2017,Lanza2019} trade convexity for a better recovery of edges and of sparse supports. The learning of a regularizer from data is now a field of its own \cite{Arridge2019,Mukherjee2023,Habring2024}. Recent contributions include convex architectures \cite{Amos2017,Goujon2023} and their weakly convex extensions \cite{Goujon2024,Shumaylov2024}, adversarial and network-based regularizers \cite{Lunz2018,Li2020NETT,Kobler2022}, learned priors and bilevel schemes \cite{Roth2009,Kunisch2013,Chen2017TNRD}, and the plug-and-play methods that replace the regularizer by a learned denoiser \cite{Venkatakrishnan2013,Romano2017,Ryu2019,Hurault2022a,Kamilov2023,Hauptmann2025}. The role of amplitude equivariance in this context, and the universality of the polyhedral norms for the convex regularizers, are established in \cite{UnserDucotterd2026}, of which the present paper is the non-convex counterpart; the gauges of star bodies also appear as the natural regularizers in the critic-based learning of \cite{Leong2024,Leong2026}. The consequences of the homogeneity of $q$ for the variational problem \eqref{Eq:VarReconstruction}---existence of a minimizer, localization of the critical set, and convergence of the steepest descent---are the subject of the companion paper \cite{Unser2026b}.

On the mathematical side, the finite-dimensional theory of the quasi-norms that Section~\ref{Sec:QuasiNorms} revisits is part of the folklore of the quasi-Banach spaces \cite{Kalton2003}, and the statements that hold on an arbitrary vector space are classical \cite{Koethe1983}; what we add is the comparison with the distance to the null space, which characterizes the quasi-seminorms within the CH class, and the structure theorem that follows from it. The star-shaped bodies and their radial functions are standard objects of convex and integral geometry \cite[Chapter 0]{Gardner2006}. The conic splines are the piecewise-linear functions on a fan of toric geometry \cite[Section 3.4]{Fulton1993} and, for $d=3$, the degree-one spherical splines of Alfeld, Neamtu and Schumaker \cite{Alfeld1996a,Alfeld1996c,Alfeld1996b}, whose approximation order is due to Neamtu and Schumaker \cite{Neamtu2004}.
The context differs on two counts: spherical splines fit data that are given on the sphere $\mathbb{S}^{2}$, a temperature map for instance, with a degree that is free, whereas we fit a quasi-norm, whose $1$-homogeneity excludes the higher degrees but allows the theory to be formulated in an arbitrary dimension $d$. What is specific to our setting is the identification of the non-negative conic splines with the polyhedral CH functionals, the reading of the nodal values as the parameters of a learning problem, the extension of the rate to the non-convex (star-shaped) bodies with a bound on the constant in terms of the Hessian of the gauge, and the second-order estimate under a measure-valued Hessian through the spherical Hessian variation, which is the restriction to the sphere of the Hessian--Schatten total variation of \cite{Aziznejad2023a}.

\subsection{Organization of the Paper}
The paper devotes one section to each of the three main issues, followed by the detailed treatment of the planar case.
In Section~\ref{Sec:QuasiNorms}, we identify the relevant class of amplitude-equivariant regularizers by linking it to the CH class (Section~\ref{Sec:AmplitudeEquivariance}) and by providing a complete characterization of the most important members of the family: the quasi-norms and the quasi-seminorms, which offer control over the null space (Section~\ref{Sec:QuasiSeminorms}). These functionals are the non-convex counterparts of the norms and the seminorms on $\R^d$, respectively.
We also make the correspondence with the star-shaped bodies and their gauges explicit (Section~\ref{Sec:StarShaped}), which is the geometric language of the approximation results. In Section~\ref{Sec:PolyhedralQuasiNorms}, we turn to the parametrization: we recall the polyhedral fans and construct the conic B-splines, with their dual basis and their stability bounds (Section~\ref{Sec:Fans}), and we prove the representation theorem for the polyhedral CH functionals (Section~\ref{Sec:PolyRepresentation}). In Section~\ref{Sec:Regularity}, we establish the rates: after defining the regularity indices of a CH functional, including the spherical Hessian variation by duality (Section~\ref{Sec:RegularityIndices}),
we prove the uniform second-order estimate together with its sharpness (Section~\ref{Sec:Approximation}) and the $L_1$ estimate under a bounded spherical Hessian variation (Section~\ref{Sec:ApproximationHTV}). In Section~\ref{Sec:Planar}, we work out the details for the planar case $d=2$, where everything becomes explicit: the conic B-splines on a uniform polar grid are trigonometric hat functions, the spherical Hessian variation of a polyhedral profile is an $\ell_1$-norm of second differences of the nodal values, and the approximation constants are sharp. 

\subsection{Notations}
\label{Sec:Notations}

The symbol $\Spc X$ denotes an abstract vector space with generic elements $x,y \in \Spc X$. The Minkowski sum of two subsets $B_1, B_2 \subseteq \Spc X$ is $B_1+B_2=\{x + y: x \in B_1, y \in B_2\}$, and the rescaling of $B \subseteq \Spc X$ by a scalar $\lambda>0$ is $\lambda B=\{\lambda x: x \in B\}$; in particular, if $B_{\|\cdot\|}=\{x \in \Spc X: \|x\|\le 1\}$ is the unit ball of a norm $\|\cdot\|$ on $\Spc X$, then $\lambda B_{\|\cdot\|}=\{x \in \Spc X: \|x\|\le \lambda\}$ is its $\lambda$-inflated version. 

When the setting is concrete with $\Spc X=\R^d$, we use boldface variables ($\V x, \V y \in \R^d$) to remind us that the objects are ``ordinary'' vectors and that the underlying results are restricted to finite dimensions. The Euclidean norm of $\R^d$ is $\|\cdot\|_2$, its unit ball is $B_{\ell_2}=\{\V x \in \R^d: \|\V x\|_2\le 1\}$, and the corresponding hypersphere is $\mathbb{S}^{d-1}=\{\V \xi \in \R^d: \|\V \xi\|_2=1\}$, so that a point $\V x\in \R^d \backslash\{\V 0\}$ is written in hyperpolar coordinates as $\V x=r \V \xi$ with $r=\|\V x\|_2 \in \R_+=[0,\infty)$ and $\V \xi=\V x/\|\V x\|_2 \in \mathbb{S}^{d-1}$. 

The {\em radial function} of a star-shaped body $B$ is denoted by
$r_B(\V \xi)$ (written $\rho_B$ in \cite{Gardner2006}): it is the radius of $B$ in the direction $\V \xi \in \mathbb{S}^{d-1}$. If $B_q$ is the unit ball of a quasi-norm $q$, then $
r_{B_q}(\V \xi)=1/q(\V \xi)$ is the reciprocal of the {\em angular profile} $q(\V \xi)$ of $q$. For $d=2$, the latter is denoted by $\varrho(\theta)=q(\V \xi(\theta))$ with $\V \xi(\theta)=(\cos \theta, \sin \theta) \in \mathbb{S}^{1}$.

The symbol $\Spc H^{k}$ denotes the $k$-dimensional Hausdorff measure---the surface measure of $\mathbb S^{d-1}$ when $k=d-1$---and is not to be confused with the Hessian $\Op H\{q\}$.

\section{Amplitude Equivariance and Quasi-norms}
\label{Sec:QuasiNorms}

\subsection{Amplitude Equivariance and the CH Class}
\label{Sec:AmplitudeEquivariance}
Following \cite{UnserDucotterd2026}, we say that a functional $g:\R^d\to[0,\infty)$ is {\em amplitude-equivariant} if there exists a function $a:\R\backslash\{0\}\to(0,\infty)$ such that
\begin{align}
\label{Eq:AEDefinition}
g(\alpha \V x)=a(\alpha)\, g(\V x), \qquad \forall \V x \in \R^d, \ \forall \alpha \in \R\backslash\{0\}.
\end{align}
The property is tailored to the variational problem \eqref{Eq:VarReconstruction}: if $g_{\V \theta}$ satisfies \eqref{Eq:AEDefinition}, then the substitution $\V s=\alpha\V s'$ shows that the reconstruction from the amplified data $\V y'=\alpha\V y$ with the adjusted regularization parameter $\lambda'=\lambda\,\alpha^2/a(\alpha)$ is $\alpha$ times the reconstruction from $\V y$ with the reference parameter $\lambda$. The functionals that we are after are those that are, in addition, continuous and symmetric. They turn out to be the powers of the members of the following class.

\begin{definition}[CH class]
\label{Def:CHscalar}
A functional $q:\R^d\to[0,\infty)$ is of class CH if it is continuous and absolutely homogeneous, in the sense that $q(\alpha \V x)=|\alpha|\,q(\V x)$ for every $\V x \in \R^d$ and every $\alpha \in \R$. We then write $q \in \CH(\R^d \to \R_+)$.
\end{definition}

A CH functional is encoded by its angular profile. Since $q(\V 0)=0$ by homogeneity, we have that
\begin{align}
\label{Eq:HomogeneousExtension}
q(\V x)= \|\V x\|_2\; q\left( \frac{\V x}{\|\V x\|_2} \right)=r \, q\left( \V \xi \right), \qquad \V x=r\V \xi \ne \V 0,
\end{align}
where $q|_{\mathbb S^{d-1}}$ is continuous and even. Conversely, the right-hand side of \eqref{Eq:HomogeneousExtension} defines a member of $\CH(\R^d\to\R_+)$ for {\em every} continuous, even, non-negative function on $\mathbb S^{d-1}$, the continuity at the origin being a consequence of the boundedness of the profile. The CH class is therefore as large as $C(\mathbb S^{d-1};\R_+)$ allows, and the only quantitative information that comes for free is an upper bound: for any norm $\|\cdot\|$ on $\R^d$, the Weierstrass theorem applied on the compact unit sphere $S_{\|\cdot\|}=\{\V x: \|\V x\|=1\}$ provides $M=\max_{\|\V x\|=1}q(\V x)<\infty$, and the homogeneity extends the bound to
\begin{equation}
\label{Eq:NormUpper}
\forall \V x \in \R^d:\quad q(\V x) \le M \|\V x\|.
\end{equation}

\begin{proposition}[Admissible amplitude-equivariant regularizers]
\label{Prop:AERegularizers}
Let $g:\R^d\to[0,\infty)$ be continuous, symmetric, amplitude-equivariant, and non-constant. Then, there exists a unique $\gamma>0$ such that
\begin{align}
\label{Eq:AEFactorization}
g=q^\gamma, \qquad q:=g^{1/\gamma},
\end{align}
where $q \in \CH(\R^d \to \R_+)$. Conversely, any functional of the form \eqref{Eq:AEFactorization} with $\gamma>0$ and $q \in \CH(\R^d \to \R_+)$ is continuous, symmetric, and amplitude-equivariant. Moreover, $g$ is convex if and only if $\gamma\ge1$ and $q$ is a seminorm.
\end{proposition}

\begin{proof}
Let us assume that $g$ satisfies the hypotheses. Since $g$ is non-constant, there exists $\V x_0$ with $g(\V x_0)>0$, so that \eqref{Eq:AEDefinition} determines $a(\alpha)=g(\alpha \V x_0)/g(\V x_0)$; this shows that $a$ is continuous and non-negative. The relation $g\big((\alpha_1\alpha_2)\V x\big)=a(\alpha_1)a(\alpha_2)g(\V x)$ gives the functional equation $a(\alpha_1\alpha_2)=a(\alpha_1)a(\alpha_2)$ for $\alpha_1,\alpha_2>0$.
We first observe that $a$ does not vanish: if $a(\alpha_1)=0$ for some $\alpha_1>0$, then $a\equiv0$ by the same relation, which would force $g$ to vanish identically. Since $a$ is continuous and positive, the function $\log a$ is continuous and additive on the multiplicative group $\R_{>0}$, hence of the form $\log a(\alpha)=\gamma \log \alpha$ for some $\gamma \in \R$. This yields $a(\alpha)=\alpha^\gamma$ for $\alpha>0$.

The continuity of $g$ at the origin then forces $\gamma>0$. Indeed, $g(\alpha \V x_0)=\alpha^\gamma g(\V x_0)$ has to remain bounded as $\alpha \downarrow 0$, which excludes $\gamma<0$; and $\gamma=0$ would give $g(\alpha\V x_0)=g(\V x_0)$ for all $\alpha>0$, so that the continuity at the origin would yield $g(\V 0)=g(\V x_0)$ and, by homogeneity again, $g\equiv g(\V 0)$, which is excluded because $g$ is non-constant. The symmetry of $g$ gives $a(-1)=1$ and hence $a(\alpha)=|\alpha|^\gamma$ for every $\alpha \ne 0$. The exponent $\gamma$ is unique because $\alpha \mapsto \alpha^\gamma$ determines it.

Let us now set $q=g^{1/\gamma}$, which is legitimate because $g\ge0$. The functional $q$ is continuous as the composition of the continuous maps $g$ and $t\mapsto t^{1/\gamma}$, and it is absolutely homogeneous:
\begin{align*}
q(\alpha \V x)=\big(|\alpha|^\gamma g(\V x)\big)^{1/\gamma}=|\alpha| \, q(\V x).
\end{align*}
Hence, $q \in \CH(\R^d \to \R_+)$, which establishes \eqref{Eq:AEFactorization}. The converse is immediate, since $q^\gamma$ inherits the continuity and the symmetry of $q$ and satisfies \eqref{Eq:AEDefinition} with $a(\alpha)=|\alpha|^\gamma$.

For the last statement, let us assume that $g=q^\gamma$ is convex. For any $\alpha \in (0,1)$, the convexity inequality applied to $g\big(\alpha\V x +(1-\alpha)\V y\big)$ with $\V y=\V 0$ gives $\alpha^\gamma g(\V x)=g(\alpha \V x)\le \alpha \, g(\V x)$, which forces $\gamma \ge 1$. Next, the sublevel set $B=\{\V x \in \R^d: g(\V x)\le 1\}=\{\V x \in \R^d: q(\V x)\le 1\}$ is convex because $g$ is convex. Given $\V x,\V y$ with $q(\V x),q(\V y)>0$, the points $\V x/q(\V x)$ and $\V y/q(\V y)$ belong to $B$ and so does their convex combination
\begin{align*}
\V z=\frac{q(\V x)}{q(\V x)+q(\V y)}\,\frac{\V x}{q(\V x)}+\frac{q(\V y)}{q(\V x)+q(\V y)}\,\frac{\V y}{q(\V y)}=\frac{\V x+\V y}{q(\V x)+q(\V y)},
\end{align*}
with $\V z \in B \iff q(\V z)\le 1$.
By homogeneity, the latter inequality is equivalent to the triangle inequality $q(\V x+\V y)\le q(\V x)+q(\V y)$. If $q(\V y)=0$, then $\V y/\epsilon\in B$ for every $\epsilon>0$ and the same argument with $\epsilon$ in place of $q(\V y)$ gives $q(\V x+\V y)\le q(\V x)+\epsilon$; the case $q(\V x)=0$ is symmetric. Thus, $q$ satisfies the triangle inequality, which makes it a seminorm. Conversely, if $q$ is a seminorm and $\gamma\ge1$, then $q^\gamma$ is convex as the composition of the convex functional $q$ with the convex, increasing map $t\mapsto t^\gamma$ on $\R_+$.
\end{proof}

Proposition~\ref{Prop:AERegularizers} is the non-convex counterpart of \cite[Theorem 2.3]{UnserDucotterd2026}, and the comparison is instructive. Convexity plays two roles there: it pins the order of homogeneity down to $\gamma\ge1$, and it reduces $q$ to a seminorm. Once it is dropped, the order of homogeneity is only constrained to be positive, and $q$ is any member of the CH class. The exponent is immaterial for the variational problem \eqref{Eq:VarReconstruction}, because $t\mapsto t^\gamma$ is increasing on $\R_+$ for every $\gamma>0$, so that we may always work with $q$ itself after a readjustment of $\lambda$. The whole point of the sections that follow is that the CH functionals are as amenable to parametrization and approximation as the seminorms are, the conic splines taking over from the polyhedral norms.

One feature of a regularizer deserves to be singled out at this stage: its null space $\Spc N_q=\{\V x \in \R^d: q(\V x)=0\}$, which is the set of signals that are not penalized at all. 
For a seminorm, the null space is automatically a linear subspace and the seminorm is constant along it. For a CH functional, the null space is only a closed symmetric cone, the set of directions $\V \xi\in\mathbb S^{d-1}$ at which the profile vanishes, and nothing forces it to be a subspace (Appendix~\ref{App:NotQuasiSeminorm} gives examples). The members of the CH class that behave like seminorms in this respect---a linear null space, and a first-order growth away from it---are the quasi-seminorms, which we now identify.

\subsection{Quasi-norms and Quasi-seminorms}
\label{Sec:QuasiSeminorms}
We start with the algebraic definition, which is meaningful on an arbitrary vector space and which is the one of the quasi-Banach literature \cite{Kalton2003}: a topological vector space is locally bounded if and only if its topology is generated by a quasi-norm, and such a space is $p$-normable for some $p\in(0,1]$ (Aoki--Rolewicz theorem \cite{Aoki1942,Rolewicz1957}).

\begin{definition}
\label{Def:QuasiNorm}
Let $\Spc X$ be a vector space. A functional $q:\Spc X \to [0,\infty)$ is called a \emph{quasi-seminorm} if the following two properties hold:
\begin{enumerate}
\item {\em Absolute homogeneity} (amplitude equivariance):
\begin{align*}
q(\alpha x)=|\alpha|\, q(x), \qquad \forall \alpha\in\R,\ \forall x\in \Spc X;
\end{align*}
\item {\em Quasi-triangle inequality}: there exists a constant $C\ge 1$ such that
\begin{align*}
q(x+y)\le C\big(q(x)+q(y)\big), \qquad \forall x,y\in \Spc X.
\end{align*}
The smallest such $C$ is called the \emph{modulus} of the quasi-seminorm.
\end{enumerate}
If, in addition, $q$ is positive-definite in the sense that $q(x)=0 \iff x=0$, then $q$ is called a \emph{quasi-norm}.
\end{definition}

A seminorm (resp., a norm) is the special case $C=1$. Since the quasi-triangle inequality is invariant under the rescaling $q \mapsto \lambda q$ with $\lambda>0$, the modulus is an intrinsic descriptor of the ``degree of non-convexity'' of $q$. On a general vector space, a quasi-norm is characterized algebraically by its unit ball, which has to be balanced, absorbing, star-shaped, and stable under addition up to a dilation (Theorem~\ref{Theo:QuasiNormsGauges} in Appendix~\ref{App:InfiniteDim}). No topology is involved there, so that no regularity can be inferred: unlike a seminorm, which is automatically continuous, a quasi-norm need not be, even on $\R^2$ (Appendix~\ref{App:Counterexample}). This is why continuity is built into the CH class, and why the quasi-(semi)norms of interest to us are the continuous ones. The first two theorems of this subsection tell us which members of the CH class they are. The first one disposes of the positive-definite case, which is the one that the parametrization of Section~\ref{Sec:PolyhedralQuasiNorms} is designed to produce: there, the quasi-triangle inequality comes for free.

\begin{theorem}[Quasi-norms are the positive-definite CH functionals]
\label{Theo:QuasiNormAuto}
Let $q \in \CH(\R^d \to \R_+)$ and let $\|\cdot\|$ be an arbitrary norm on $\R^d$. Then, the following statements are equivalent.
\begin{enumerate}
\item {\em Trivial null space:} $q$ is positive-definite, in the sense that $\Spc N_q=\{\V 0\}$.
\item {\em Quasi-norm:} $q$ is a continuous quasi-norm on $\R^d$.
\item {\em Two-sided norm equivalence:} there exist two constants $0<m\le M<\infty$ such that
\begin{align}
\label{Eq:TwoSided}
\forall \V x\in \R^d:\quad m\|\V x\| \le q(\V x) \le M\|\V x\|.
\end{align}
\item {\em Bounded unit ball:} the set $B_q=\{\V x \in \R^d: q(\V x)\le1\}$ is bounded.
\end{enumerate}
In that case, the optimal constants of \eqref{Eq:TwoSided} are $m=\min_{\|\V x\|=1}q(\V x)$ and $M=\max_{\|\V x\|=1}q(\V x)$; the lower bound
\begin{align}
\label{Eq:qlowerbound}
\forall \V x \in \R^d:\quad m \|\V x\| \le q(\V x)
\end{align}
is equivalent to the inclusion $B_q \subseteq \frac{1}{m} B_{\|\cdot\|}$; the unit ball satisfies $\frac{1}{M}B_{\|\cdot\|}\subseteq B_q \subseteq \frac{1}{m}B_{\|\cdot\|}$; and the modulus of $q$ admits the explicit bound
\begin{align}
\label{Eq:ModulusBound}
C\le \frac{M}{m}=\frac{\max_{\|\V x\|=1}q(\V x)}{\min_{\|\V x\|=1}q(\V x)}.
\end{align}
This implies that $q$ is topologically equivalent to any norm on $\R^d$.
\end{theorem}

\begin{proof}
The upper bound of \eqref{Eq:TwoSided} is \eqref{Eq:NormUpper}, which holds for every CH functional, and its geometric transcription is $\frac{1}{M}B_{\|\cdot\|}\subseteq B_q$. We prove $(1)\implies(3)\implies(2)\implies(1)$, and then $(3)\iff(4)$.

\smallskip
$(1) \implies (3)$: A positive-definite $q$ is continuous and strictly positive on the compact unit sphere $S_{\|\cdot\|}=\{\V x: \|\V x\|=1\}$, so that $m=\min_{\|\V x\|=1}q(\V x)>0$. Absolute homogeneity then extends the inequality $m\|\V x\|\le q(\V x)$ to the whole of $\R^d$. Both constants are optimal, again by homogeneity.

\smallskip
$(3) \implies (2)$: The triangle inequality of the norm gives
\begin{align*}
q(\V x+\V y)\le M\|\V x+\V y\|\le M\big(\|\V x\|+\|\V y\|\big)\le \frac{M}{m}\big(q(\V x)+q(\V y)\big),
\end{align*}
which is the quasi-triangle inequality with $C\le M/m$ and establishes \eqref{Eq:ModulusBound}; the lower bound forces $\Spc N_q=\{\V 0\}$, and $q$ is continuous by assumption.

\smallskip
$(2) \implies (1)$: A quasi-norm is positive-definite by Definition~\ref{Def:QuasiNorm}.

\smallskip
$(3) \iff (4)$: If \eqref{Eq:qlowerbound} holds, then every $\V x\in B_q$ satisfies $\|\V x\|\le\frac1m$, so that $B_q\subseteq\frac1mB_{\|\cdot\|}$ is bounded. Conversely, if $B_q\subseteq\frac1mB_{\|\cdot\|}$ for some $m>0$, then no $\V x\ne\V 0$ has $q(\V x)=0$, since the whole line $\R\V x$ would otherwise sit inside the bounded set $B_q$; and, for $q(\V x)>0$, the point $\V x/q(\V x)\in B_q$ satisfies $\|\V x\|/q(\V x)\le\frac1m$, which is \eqref{Eq:qlowerbound}.
\end{proof}

A positive-definite CH functional is therefore a quasi-norm by construction, with a modulus that is bounded by the ratio of the largest to the smallest value of its profile. The prototypical example is the $\ell_p$ quasi-norm $\|\V x\|_p=\big(\sum_n|x_n|^p\big)^{1/p}$ with $0<p<1$, which fails the (strict) triangle inequality as soon as $d\ge2$ (take $\V x=\V e_1$ and $\V y=\V e_2$, for which $\|\V x+\V y\|_p=2^{1/p}>2$) and which is not Lipschitz-continuous, a point to which we return in Appendix~\ref{App:ellp}. 

The next question is what happens when $q$ is allowed to vanish outside the origin, which is the configuration of the regularizers with a null space. Theorem~\ref{Theo:QuasiNormAuto} no longer applies, and the quasi-triangle inequality is no longer automatic: Appendix~\ref{App:NotQuasiSeminorm} exhibits two continuous, absolutely homogeneous, and non-negative functionals on $\R^2$ that are not quasi-seminorms---a conic spline whose null space is the union of the two diagonals, hence not a linear subspace, and a functional whose null space is a line but which is not comparable to the distance to that line. The next theorem shows that these are the only two obstructions and, more positively, that the quasi-seminorms are exactly the members of the CH class that are comparable to a seminorm: the null space is then a linear subspace and $q$ grows to the first order away from it, as a seminorm does. We recall that the distance from $\V x\in\R^d$ to a non-empty closed set $\Spc N \subseteq \R^d$ is $\dist(\V x,\Spc N)=\min_{\V y \in \Spc N}\|\V x-\V y\|$, and that this distance is a seminorm whenever $\Spc N$ is a linear subspace.

\begin{theorem}[Quasi-seminorms within the CH class]
\label{Theo:SemiQuasiNorm}
Let $q \in \CH(\R^d \to \R_+)$, let $\Spc N_q = \{\V x \in \R^d : q(\V x) = 0\}$ be its null space, and let $\|\cdot\|$ be an arbitrary norm on $\R^d$. Then, the following conditions are equivalent.
\begin{enumerate}
\item {\em Quasi-triangle inequality:} there exists a constant $C\geq 1$ such that
\begin{equation}
\label{Eq:QuasiTriangle}
\forall \V x,\V y \in \R^d: \quad q(\V x+\V y) \leq C\big(q(\V x)+q(\V y)\big).
\end{equation}
\item {\em Loosened convexity:} there exists a constant $C\geq 1$ such that
\begin{equation}
\label{Eq:QuasiConvex}
\forall \V x,\V y \in \R^d, \ \forall t\in[0,1]: \quad q\big((1-t)\V x+t\V y\big) \leq C\big((1-t)q(\V x)+tq(\V y)\big).
\end{equation}
\item {\em Structure of the null space:} the set $\Spc N_q$ is a linear (possibly trivial) subspace of $\R^d$, and there exist two constants $0<m_0\le M_0<\infty$ such that
\begin{equation}
\label{Eq:DistanceComparison}
\forall \V x \in \R^d:\quad m_0\, \dist(\V x,\Spc N_q) \le q(\V x) \le M_0\, \dist(\V x,\Spc N_q),
\end{equation}
where the metric in $\dist$ is any chosen norm $\|\cdot\|$ on $\R^d$.
\item {\em Geometry of the unit ball:} the unit ball $B_q=\{\V x\in \R^d : q(\V x) \le 1\}$ is a symmetric, closed, star-shaped set that contains $\V 0$ as an interior point and that satisfies the inclusion property $B_q + B_q \subseteq a_0 B_q$ for some $a_0\ge 2$.
\item {\em Comparability to a seminorm:} there exist a seminorm $p$ on $\R^d$ and a constant $c\ge1$ such that
\begin{equation}
\label{Eq:SeminormComparison}
\forall \V x \in \R^d:\quad \tfrac1c\,p(\V x)\le q(\V x)\le c\,p(\V x).
\end{equation}
\end{enumerate}
The functional $q$ is then a continuous quasi-seminorm on $\R^d$, the optimal constants of Items 1 and 4 being related by $C\le a_0\le 2C$. Finally, $q$ is a seminorm if and only if \eqref{Eq:QuasiTriangle} (resp., \eqref{Eq:QuasiConvex}, \eqref{Eq:SeminormComparison}) holds with $C=1$ (resp., $c=1$) or, equivalently, if and only if $B_q$ in Item 4 is convex.
\end{theorem}

\begin{proof}
We first record what holds for {\em every} CH functional: its unit ball $B_q$ is symmetric, because $q(-\V x)=q(\V x)$; closed, as the preimage $q^{-1}([0,1])$ of a closed set under a continuous map; star-shaped, because $q(t\V x)=t\,q(\V x)\le q(\V x)$ for $t\in[0,1]$; and it contains $\frac1MB_{\|\cdot\|}$ by \eqref{Eq:NormUpper}, so that $\V 0$ is an interior point. Item 4 therefore reduces to the inclusion property. We prove the cycle $(1)\implies(2)\implies(4)\implies(1)$ and then $(1)\implies(3)\implies(5)\implies(1)$.

\smallskip
$(1) \implies (2)$: By absolute homogeneity, $q\big((1-t)\V x+t\V y\big)\le C\big(q((1-t)\V x)+q(t\V y)\big)=C\big((1-t)q(\V x)+tq(\V y)\big)$, with the same constant $C$.

\smallskip
$(2) \implies (4)$: For $\V x,\V y\in B_q$, homogeneity and \eqref{Eq:QuasiConvex} with $t=\frac12$ give $q(\V x+\V y)=2\,q\big(\tfrac12\V x+\tfrac12\V y\big)\le 2C$, so that $\V x+\V y\in 2C\,B_q$: the inclusion property holds with $a_0=2C$.

\smallskip
$(4) \implies (1)$: Let $\V x,\V y\in\R^d$ and $t>0$ with $t\ge\max\{q(\V x),q(\V y)\}$. Then, $\V x/t$ and $\V y/t$ belong to $B_q$, so that $(\V x+\V y)/t\in a_0B_q$, which reads $q(\V x+\V y)\le a_0t$. The infimum over the admissible $t$ yields $q(\V x+\V y)\le a_0\max\{q(\V x),q(\V y)\}\le a_0\big(q(\V x)+q(\V y)\big)$, which is \eqref{Eq:QuasiTriangle} with $C=a_0$. Together with the value $a_0=2C$ of the previous step, this gives $C\le a_0\le2C$ for the optimal constants.

\smallskip
$(1) \implies (3)$: The set $\Spc N_q$ is closed under scalar multiplication by homogeneity and under addition by \eqref{Eq:QuasiTriangle}, hence it is a linear subspace. Let $\M P_\perp$ be the orthogonal projector onto $\Spc N_q^\perp$, so that $\dist(\V x,\Spc N_q)=\|\M P_\perp\V x\|_2$ for the Euclidean norm; the case of another norm follows from the equivalence of the norms on $\R^d$. Since $\V x-\M P_\perp\V x\in\Spc N_q$, two applications of \eqref{Eq:QuasiTriangle}, to $\V x=\M P_\perp\V x+(\V x-\M P_\perp\V x)$ and to $\M P_\perp\V x=\V x-(\V x-\M P_\perp\V x)$, give
\begin{align}
\label{Eq:ProjectionComparison}
\forall \V x\in\R^d:\quad \tfrac1C\,q(\M P_\perp\V x)\ \le\ q(\V x)\ \le\ C\,q(\M P_\perp\V x).
\end{align}
The restriction $\bar q=q|_{\Spc N_q^\perp}$ is a CH functional on $\Spc N_q^\perp$ that is positive-definite there, since $\Spc N_q\cap\Spc N_q^\perp=\{\V 0\}$; Theorem~\ref{Theo:QuasiNormAuto}, applied on the space $\Spc N_q^\perp$, gives $m\|\V z\|_2\le q(\V z)\le M\|\V z\|_2$ for every $\V z\in\Spc N_q^\perp$. With $\V z=\M P_\perp\V x$, \eqref{Eq:ProjectionComparison} yields \eqref{Eq:DistanceComparison} with $m_0=m/C$ and $M_0=CM$.

\smallskip
$(3) \implies (5)$: Since $\Spc N_q$ is a linear subspace, $p=\dist(\cdot,\Spc N_q)$ is a seminorm, and \eqref{Eq:DistanceComparison} is \eqref{Eq:SeminormComparison} with $c=\max\{M_0,1/m_0\}$.

\smallskip
$(5) \implies (1)$: The triangle inequality of $p$ gives $q(\V x+\V y)\le c\,p(\V x+\V y)\le c\big(p(\V x)+p(\V y)\big)\le c^2\big(q(\V x)+q(\V y)\big)$, which is \eqref{Eq:QuasiTriangle} with $C=c^2$.

\smallskip
{\em The seminorm case}: For $C=1$, \eqref{Eq:QuasiTriangle} is the triangle inequality and $q$ is a seminorm (and $c=1$ in \eqref{Eq:SeminormComparison} with $p=q$), whose unit ball is convex as a sublevel set of a convex function. Conversely, if $B_q$ is convex, then the argument of the proof of Proposition~\ref{Prop:AERegularizers}---the convex combination of $\V x/q(\V x)$ and $\V y/q(\V y)$---shows that $q$ satisfies the triangle inequality.
\end{proof}

Theorem~\ref{Theo:QuasiNormAuto} is the case $\Spc N_q=\{\V 0\}$ of Theorem~\ref{Theo:SemiQuasiNorm}, Item 3, for which \eqref{Eq:DistanceComparison} is \eqref{Eq:TwoSided}. In the convex theory, a seminorm is {\em exactly} constant along its null space, so that it is the composition of a norm on the quotient $\R^d/\Spc N_q$ with the quotient map. The two-sided inequality \eqref{Eq:ProjectionComparison} is the non-convex substitute for that quotient argument, and it only holds up to the factor $C$. The quasi-seminorms for which it holds with equality are the ones that we recommend for the design of a regularizer, and we give them a name.

\begin{definition}[Null-space control]
\label{Def:NullSpaceControl}
A functional $q:\R^d\to[0,\infty)$ is a {\em quasi-seminorm with null-space control} if there exist a matrix $\M P\in\R^{d'\times d}$ and a continuous quasi-norm $\bar q$ on $\R^{d'}$ such that
\begin{align}
\label{Eq:NullSpaceControl}
q=\bar q\circ\M P;\qquad\mbox{that is,}\qquad q(\V x)=\bar q(\M P\V x), \quad \forall \V x \in \R^d.
\end{align}
\end{definition}

Such a functional is a continuous quasi-seminorm with $\Spc N_q=\ker(\M P)$, by the linearity of $\M P$ and the quasi-triangle inequality of $\bar q$; its modulus does not exceed that of $\bar q$. The null space is {\em controlled} in the sense that it is set by $\M P$ alone---a finite-difference operator or a filter bank, say, that encodes the invariances one wants to impose---while $\bar q$ is positive-definite and therefore, by Theorem~\ref{Theo:QuasiNormAuto}, admissible without any further condition. This is the exact analogue of the analysis form of the convex theory \cite{UnserDucotterd2026}, in which a norm is composed with a trainable regularization operator. The next theorem characterizes the quasi-seminorms with null-space control within the CH class and shows that every other quasi-seminorm is comparable to one of them.

\begin{theorem}[Structure of continuous quasi-seminorms]
\label{Theo:QuasiSeminormStructure}
Let $q\in\CH(\R^d\to\R_+)$ and let $\Spc N_q$ be its null space.
\begin{enumerate}
\item {\em Null-space control.} The following three conditions are equivalent:
\begin{enumerate}
\item[(i)] $q$ is a quasi-seminorm with null-space control, as in \eqref{Eq:NullSpaceControl};
\item[(ii)] $q$ is invariant along its null space, in the sense that $q(\V x+\V y)=q(\V x)$ for every $\V x \in \R^d$ and every $\V y \in \Spc N_q$;
\item[(iii)] $q=\bar q\circ\M P_\perp$, where $\Spc N_q$ is a linear subspace, $\M P_\perp$ is the orthogonal projector onto $\Spc N_q^\perp$, and $\bar q=q|_{\Spc N_q^\perp}$ is a continuous quasi-norm on $\Spc N_q^\perp$.
\end{enumerate}
In that case, $\Spc N_q=\ker(\M P)$, the modulus of $q$ does not exceed that of $\bar q$, and the unit ball is the cylinder $B_q=B_{\bar q}+\Spc N_q$, where $B_{\bar q}$ is the unit ball of $q|_{\Spc N_q^\perp}$ in $\Spc N_q^\perp$.
\item {\em General case.} The functional $q$ is a quasi-seminorm if and only if there exist a matrix $\M P\in\R^{d'\times d}$, a continuous quasi-norm $\bar q$ on $\R^{d'}$, and a constant $c\ge1$ such that
\begin{align}
\label{Eq:QuotientComparison}
\forall \V x \in \R^d: \quad \tfrac1c\,\bar q(\M P \V x) \le q(\V x) \le c\,\bar q(\M P \V x).
\end{align}
In that case, $\Spc N_q=\ker(\M P)$. If $q$ is a quasi-seminorm with modulus $C$, then \eqref{Eq:QuotientComparison} holds with $c=C$, $\M P=\M P_\perp$, and $\bar q=q|_{\Spc N_q^\perp}$; conversely, if $\bar q$ has modulus $\bar C$ in \eqref{Eq:QuotientComparison}, then the modulus of $q$ does not exceed $c^2\bar C$.
\end{enumerate}
\end{theorem}

\begin{proof}
{\em Item 1}: If $q=\bar q\circ\M P$, then $\Spc N_q=\ker(\M P)$ because $\bar q$ is positive-definite, and $q(\V x+\V n)=\bar q(\M P\V x+\M P\V n)=q(\V x)$ for every $\V n\in\ker(\M P)$, which is (i)$\implies$(ii). If $q$ is invariant along its null space, then $\Spc N_q$ is closed under scalar multiplication by homogeneity and under addition by invariance, hence it is a linear subspace; the restriction $\bar q=q|_{\Spc N_q^\perp}$ is a continuous quasi-norm on $\Spc N_q^\perp$, exactly as in the proof of Theorem~\ref{Theo:SemiQuasiNorm}, $(1)\implies(3)$; and the invariance gives $q(\V x)=q\big(\M P_\perp\V x+(\V x-\M P_\perp\V x)\big)=\bar q(\M P_\perp\V x)$, which is (ii)$\implies$(iii). The implication (iii)$\implies$(i) is immediate. In that case, $q(\V x+\V y)=\bar q(\M P\V x+\M P\V y)\le\bar C\big(q(\V x)+q(\V y)\big)$, where $\bar C$ is the modulus of $\bar q$, shows that the modulus of $q$ does not exceed that of $\bar q$, and $q(\V x)\le1$ holds if and only if $\M P_\perp\V x\in B_{\bar q}$, that is, if and only if $\V x\in B_{\bar q}+\Spc N_q$.

\smallskip
{\em Item 2}: If $q$ is a quasi-seminorm with modulus $C$, then \eqref{Eq:ProjectionComparison} is \eqref{Eq:QuotientComparison} with $c=C$, $\M P=\M P_\perp$, and $\bar q=q|_{\Spc N_q^\perp}$, which is a continuous quasi-norm by the proof of Theorem~\ref{Theo:SemiQuasiNorm}. Conversely, let \eqref{Eq:QuotientComparison} hold. Since $\bar q$ is positive-definite, $q(\V x)=0$ if and only if $\M P\V x=\V 0$, so that $\Spc N_q=\ker(\M P)$ is a linear subspace, and the linearity of $\M P$ with the quasi-triangle inequality of $\bar q$ gives
\begin{align*}
q(\V x+\V y)\le c\,\bar q(\M P\V x+\M P\V y)\le c\bar C\big(\bar q(\M P\V x)+\bar q(\M P\V y)\big)\le c^2\bar C\big(q(\V x)+q(\V y)\big),
\end{align*}
which is \eqref{Eq:QuasiTriangle} with $C\le c^2\bar C$.
\end{proof}

The quasi-seminorms that are not of the controlled type are somewhat exotic. For instance, the profile $\varrho(\theta)=|\sin\theta|\,(2+\cos2\theta)$ is even and $\pi$-periodic and vanishes to the first order at $\theta=0$, so that its homogeneous extension $q(r\V \xi(\theta))=r\varrho(\theta)$ is a quasi-seminorm on $\R^2$ with $\Spc N_q=\R\V e_1$ and a modulus that does not exceed $3$; yet $q(0,1)=1$ and $q(3,1)=2.8$, although the two points are at the same distance from $\Spc N_q$. Such a functional penalizes a signal differently according to its component along the null space, which defeats the purpose of having one. For the design of a regularizer, the practical consequence is that one should not attempt to create a null space by letting a non-convex functional degenerate: the admissibility is then at the mercy of Theorem~\ref{Theo:SemiQuasiNorm}, Item 3, and the invariance along the null space is lost in general. The safe construction is \eqref{Eq:NullSpaceControl}, where the matrix $\M P$ carries the null space and where the conic splines of Section~\ref{Sec:PolyhedralQuasiNorms} parametrize the positive-definite factor $\bar q$, so that no condition beyond the positivity of the nodal values has to be enforced. The consequences of these results for the variational problem \eqref{Eq:VarReconstruction}---the existence of a minimizer and the convergence of descent algorithms to a critical point---follow from the homogeneity of $q$ alone, without any convexity; they are developed in the companion paper \cite{Unser2026b}.

\subsection{Unit Balls, Radial Functions, and Gauges}
\label{Sec:StarShaped}
Theorem~\ref{Theo:SemiQuasiNorm} shows that, as is the case for norms, quasi-norms are characterized by their unit ball. We now make the correspondence explicit, since it is the geometric side of the approximation theory of Section~\ref{Sec:Regularity}: the conic spline interpolant of a convex $q$ is the gauge of an inscribed polytope, and the general case is the extension of that picture to star-shaped bodies.

A subset $S \subseteq \R^d$ is \emph{star-shaped} (with respect to the origin) if $t\V x\in S$ for every $\V x\in S$ and every $t\in[0,1]$. The Banach disk $B_{\|\cdot\|}=\{\V x \in \R^d: \|\V x\|\le 1\}$ of a norm is the classic example; it is star-shaped not only with respect to the origin but with respect to each of its points, which is the same as being convex. The star-shaped bodies of interest to us are the unit balls $B_q$ of the continuous quasi-norms $q$. The quasi-sphere associated with $q$ is the set $\partial B_q=\{\V x \in \R^d: q(\V x)=1\}$, which is encoded by the radial function
\begin{align}
\label{Eq:RadialProfile}
r_{B_q}(\V \xi)=\sup\{r\ge 0:\ r\V \xi\in B_q\}=
\begin{cases}\frac{1}{q(\V \xi)}, & q(\V \xi)>0,\\
+\infty,& \mbox{otherwise,}\end{cases}
\end{align}
parametrized by the unit vectors $\V \xi \in \mathbb{S}^{d-1}$. The second branch of \eqref{Eq:RadialProfile} occurs only for a quasi-seminorm with a non-trivial null space, in which case the unit ball is unbounded in the directions of $\Spc N_q$; for a quasi-norm, the positive-definiteness of $q$ ensures that $q(\V \xi)$ does not vanish and the radial function is finite. The homogeneous extension \eqref{Eq:HomogeneousExtension} then reads $q(\V x)=r/r_{B_q}(\V \xi)$ for $\V x=r\V \xi$: the angular profile of $q$ and the radial function of its unit ball are reciprocal to each other.

The converse direction goes through the \emph{gauge} (or Minkowski functional) of a star-shaped set $S\subseteq \R^d$ that contains a neighborhood of the origin,
\begin{align}
\label{Eq:Gauge}
\Gamma_S : \R^d \to [0,+\infty), \qquad \Gamma_S(\V x) \coloneqq \inf\{ \lambda > 0 : \V x \in \lambda S \}.
\end{align}
The next theorem is the natural extension of a classical result in Banach space theory; see Theorem~\ref{Theo:QuasiNormsGauges} in Appendix~\ref{App:InfiniteDim} for the more general---but weaker---statement that holds on an arbitrary vector space. The crucial ingredient is again the continuity of $q$, which brings a refinement over the description of $B_q$ in Theorem~\ref{Theo:SemiQuasiNorm}, Item 4: the radial function is continuous, and the inclusion property $B_q+B_q\subseteq a_0B_q$ becomes a consequence rather than a hypothesis.

\begin{theorem}[Gauge--quasi-norm equivalence]
\label{Theo:GaugeSemiQuasiNorm}
Let $S\subseteq \R^d$. The following statements are equivalent.
\begin{enumerate}
\item There exists a continuous quasi-norm $q:\R^d\to[0,\infty)$ such that
\begin{align*}
S=\{\V x \in \R^d:\ q(\V x)\le 1\}.
\end{align*}
\item $S$ is a symmetric, closed, bounded,
star-shaped body whose radial function is continuous and strictly positive.
\end{enumerate}
In that case, $q$ is uniquely given by the gauge $q=\Gamma_S$ specified by \eqref{Eq:Gauge}.
\end{theorem}

\begin{proof}
The proof is split into three parts.

\medskip
\emph{Step 1, $(1) \implies (2)$}:
Let us assume first that $S=\{\V x:\ q(\V x)\le 1\}$, where $q$ is a continuous quasi-norm. Then, $S$ is closed by continuity of $q$, symmetric because $q(-\V x)=q(\V x)$, and star-shaped because $q(t\V x)=t \, q(\V x)\le q(\V x)$ for $0\le t\le 1$. Its boundedness and its inclusion of a neighborhood of the origin follow from the relation $\frac{1}{M}B_{\|\cdot\|}\subseteq S\subseteq \frac{1}{m}B_{\|\cdot\|}$ of Theorem~\ref{Theo:QuasiNormAuto}.
The radial function of $S$, see \eqref{Eq:RadialProfile}, is given by
\begin{align*}
r_S(\V \xi)=\frac{1}{q(\V \xi)},\quad \V \xi \in \mathbb S^{d-1},
\end{align*}
and it is continuous and strictly positive, since $q$ is continuous and does not vanish on $\mathbb S^{d-1}$,
with $r_S\ge \frac1M$. For later use, we also record the inclusion property: if $q(\V x+\V y)\le C\big(q(\V x)+q(\V y)\big)$, then $q(\V x+\V y)\le 2C$ for $\V x,\V y\in S$, so that $S+S\subseteq 2C S$ with $2C\ge2$.

\medskip
\emph{Step 2, $(2) \implies (1)$}:
Under the assumption that $S$ satisfies the conditions in (2), the radial function $r_S$ attains a positive minimum $r_{\min}$ on the compact sphere, so that $r_{\min}B_{\ell_2}\subseteq S$ by star-shapedness: $S$ contains a neighborhood of the origin and the gauge $q = \Gamma_S$ is well-defined.

\smallskip
\noindent\textit{A. Homogeneity and symmetry}. Positive homogeneity $q(\alpha \V x) = \alpha \, q(\V x)$ for $\alpha > 0$ follows directly from the scaling of the set $\lambda S$, while the symmetry $q(-\V x) = q(\V x)$ results from $S = -S$. Thus, $q(\alpha \V x) = |\alpha| \, q(\V x)$ for all $\alpha \in \R$.

\smallskip
\noindent\textit{B. Continuity}.
Since $S$ is star-shaped with a continuous radial function $r_S$, we have, for $\V x\neq \V 0$,
\begin{align*}
\Gamma_S(\V x) = \frac{\|\V x\|_2}{r_S(\V x/\|\V x\|_2)}, \qquad \Gamma_S(\V 0)=0.
\end{align*}
Because $r_S$ is continuous and strictly positive on the compact sphere, the continuity of $r_S$ is transferred to $\Gamma_S$, first on $\mathbb{S}^{d-1}$ and then on $\R^d \backslash \{\V 0\}$ by homogeneous extension. At the origin, the inclusion $r_{\min}B_{\ell_2}\subseteq S$ gives $\epsilon r_{\min}B_{\ell_2}\subseteq\epsilon S$ for every $\epsilon>0$, so that $\|\V x\|_2<\epsilon r_{\min}$ implies $\Gamma_S(\V x)\le\epsilon$, which proves the continuity of $\Gamma_S$ at $\V x=\V 0$.

\smallskip
\noindent\textit{C. Quasi-triangle inequality}. The boundedness of $S$ gives $S\subseteq \frac{1}{m}B_{\|\cdot\|}$ for some $m>0$, which yields the lower bound \eqref{Eq:qlowerbound} and hence the positive-definiteness of $q$. Theorem~\ref{Theo:QuasiNormAuto} then tells us that $q$ is a continuous quasi-norm, with a modulus that does not exceed $M/m$. In particular, the inclusion property $S+S\subseteq 2C S$ of Theorem~\ref{Theo:SemiQuasiNorm}, Item 4, is a consequence of the hypotheses of Item 2 of the present theorem rather than an additional requirement.

\medskip
\emph{Step 3, set reconstruction and uniqueness}:
We first prove that $S = \{\V x\in \R^d : \Gamma_S(\V x) \le 1\}$ by establishing a double inclusion. If $\V x\in S$, then $\V x\in 1 \cdot S$, which, by the definition of $\Gamma_S$, implies that $S \subseteq \{\V x \in \R^d: \Gamma_S(\V x) \le 1\}$. Conversely, let us take any fixed $\V x$ such that $\Gamma_S(\V x) \le 1$. By the properties of an infimum, there exists a decreasing sequence of positive scalars $t_n \to \Gamma_S(\V x)$ such that $\V x \in t_n S$ or, equivalently, $\frac{\V x}{t_n} \in S$.
\begin{itemize}
    \item If $\Gamma_S(\V x) < 1$, then $t_n < 1$ for $n$ sufficiently large. Because $S$ is star-shaped and $\V 0 \in S$, the scaled vector $\V x = t_n \left(\frac{\V x}{t_n}\right)$ belongs to $S$.
    \item If $\Gamma_S(\V x) = 1$, then the sequence of vectors $\V y_n = \frac{\V x}{t_n} \in S$ converges to $\V x$ as $t_n \to 1^+$. Because $S$ is closed, it contains its limit points, which forces $\V x = \lim_{n \to \infty} \V y_n \in S$.
\end{itemize}
Thus, $\{ \V x \in \R^d : \Gamma_S(\V x) \le 1 \} \subseteq S$ and, consequently, $S=\{\V x \in \R^d: q(\V x)\le 1\}$. By invoking the homogeneity of $q$, this gives
\begin{align*}
\Gamma_S(\V x) = \inf\{\lambda>0:\ q(\V x)\le \lambda\} = q(\V x),
\end{align*}
which proves uniqueness.
\end{proof}

In the language of \cite{Gardner2006}, Theorem~\ref{Theo:GaugeSemiQuasiNorm} states that the continuous quasi-norms on $\R^d$ are the gauges of the symmetric star bodies with a continuous positive radial function, and Theorem~\ref{Theo:SemiQuasiNorm} adds that the continuous quasi-seminorms are, among the gauges of the symmetric closed star-shaped sets that contain the origin in their interior (bounded or not), those whose set satisfies $B+B\subseteq a_0B$. The polyhedral members of the class are the gauges of the symmetric star-shaped polytopes, which is the geometric form of the representation theorem of the next section (Theorem~\ref{Theo:PolyQuasiSemiNorms}, Item 4).

\section{Conic Splines and Polyhedral Quasi-norms}
\label{Sec:PolyhedralQuasiNorms}

\subsection{Polyhedral Fans and Conic B-splines}
\label{Sec:Fans}

A \emph{polyhedral cone} in $\R^d$ is a set of the form
\begin{align}
C=\Cone\{\M r_1,\ldots,\M r_k\} = \left\{ \sum_{i=1}^k \lambda_i \M r_i: \lambda_i\ge 0 \right\},
\end{align}
where $\M r_1,\ldots,\M r_k\in\R^d$. Such a cone is called \emph{simplicial} if its generators $\M r_1,\ldots,\M r_k$ are linearly independent, in which case its dimension is $\dim(C)=k$.

\begin{definition}[Polyhedral fan]
\label{Def:Fan}
A \emph{polyhedral fan} in $\R^d$ is a finite collection $\Spc F$ of polyhedral cones such that
\begin{enumerate}
\item every face of a cone in $\Spc F$ also belongs to $\Spc F$;
\item for any $C_1,C_2\in\Spc F$, the intersection $C_1\cap C_2$ is a common face of both cones.
\end{enumerate}
The fan is said to be
\begin{itemize}
\item \emph{complete} if its cones cover the whole space with $\bigcup_{C_{\alpha}\in\Spc F} C_{\alpha}=\R^d$;
\item \emph{simplicial} if all its cones are simplicial.
\end{itemize}
\end{definition}

A standard result in polyhedral geometry is that every polyhedral fan admits a simplicial refinement; see Proposition~\ref{Prop:FanSimplicial} in Appendix~\ref{App:Polyhedral}. There is therefore no loss of generality in focusing on complete simplicial fans.

A cone $C\in\Spc F$ is \emph{maximal} if there is no other cone $\widetilde C\in\Spc F \backslash\{C\}$ such that $C\subseteq \widetilde C$. In a complete $d$-dimensional simplicial fan, the maximal cones are exactly the $d$-dimensional cones; that is, the {\em cells} of $\Spc F$. The cones of dimension $(d-1)$ are called {\em facets}. Since each facet in $\Spc F$ is uniquely identified as the intersection of two cells, we shall generally denote it by $F_{\alpha\beta}=C_{\alpha}\cap C_{\beta}$.

Let $\Spc F$ be a complete simplicial fan in $\R^d$ with normalized rays $\{\M r_1, \ldots, \M r_N\} \subset \mathbb{S}^{d-1}$. Each $d$-dimensional simplicial cell $C_\alpha \in \Spc F$ is the conic hull of its generating rays,
\begin{align}
C_\alpha = \Cone\{\M r_{\alpha,1}, \ldots, \M r_{\alpha,d}\},
\end{align}
and is described by its vertex matrix
\begin{align}
\label{Eq:VertexMatrix}
\M R_\alpha = \begin{bmatrix} \M r_{\alpha,1} & \cdots & \M r_{\alpha,d} \end{bmatrix} \in \R^{d \times d}.
\end{align}

The space of continuous functions on $\R^d$ that are piecewise-linear over the cells of the fan $\Spc F$ is denoted by $\CPWL_{\Spc F}$ and is defined as
\begin{align}
\label{Eq:CPWLF}
\CPWL_{\Spc F} = \left\{ f \in C(\R^d) : \forall C_\alpha \in \Spc F, \ f|_{C_\alpha}(\V x) = \langle\M f_\alpha, \V x\rangle \text{ for some } \M f_\alpha \in \R^d \right\}.
\end{align}
We shall refer to the elements of this space as \emph{conic splines}.

The conic B-spline associated with the ray $\M r_i$ is the unique function $\varphi_i \in \CPWL_{\Spc F}$ such that
\begin{align}
\label{Eq:ConicBsplineNodal}
\varphi_i(\M r_j)=\delta_{ij}.
\end{align}
These functions are of special interest because they provide a representation of conic splines that is stable in the sense of Riesz,
as stated in Proposition~\ref{Prop:ConicBsplines}, Item 4. They are the conic analogues of the Courant hat functions that are widely used in finite-element analysis.

\begin{proposition}[Conic B-spline representation]
\label{Prop:ConicBsplines}
Let $\Spc F$ be a complete simplicial fan in $\R^d$ with normalized rays $\{\M r_i\}_{i=1}^N$. To each cell $C_\alpha = \Cone\{\M r_{\alpha,1},\ldots,\M r_{\alpha,d}\} \in \Spc F$, we associate a dual system $\{\tilde{\M r}_{\alpha,1},\dots,\tilde{\M r}_{\alpha,d}\}\subset \R^d$ that is uniquely specified by the biorthogonality condition
\begin{align}
\label{Eq:Biorthogonality}
\tilde{\M r}_{\alpha,j}^{\Tr}\M r_{\alpha,k} = \delta_{jk}.
\end{align}
Then, the following holds.
\begin{enumerate}
\item {\em Explicit form}: the conic B-spline associated with the $i$th ray of $\Spc F$ is
\begin{align}
\label{Eq:ConicBspline}
\varphi_i(\V x)=\begin{cases}
\tilde{\M r}_{\alpha,k}^{\Tr}\V x, & \text{if } \V x\in C_\alpha \text{ and } \M r_i=\M r_{\alpha,k},
\\[1ex]
0, & \text{if } \M r_i\notin C_\alpha,
\end{cases}
\end{align}
with its support being the closed star $\operatorname{supp}(\varphi_i) =\St(\M r_i):= \bigcup_{\alpha:\,\M r_i\in C_\alpha} C_\alpha$ of $\M r_i$ in $\Spc F$. It is non-negative and positively homogeneous,
\begin{align}
\label{Eq:BsplinePositivity}
\varphi_i(\V x)&\ge 0,\\
\label{Eq:BsplineHomogeneity}
\varphi_i(\lambda \V x) &= \lambda \varphi_i(\V x),
\end{align}
for all $\V x \in \R^d$ and $\lambda \ge 0$.
\item {\em Nodal basis}: every conic spline $s \in \CPWL_{\Spc F}$ admits the expansion
\begin{align}
\label{Eq:NodalExpansion}
s(\V x) = \sum_{i=1}^{N} s(\M r_i)\,\varphi_i(\V x)
\end{align}
and is uniquely determined by its nodal values $\{s(\M r_1),\dots,s(\M r_N)\}$.
\item {\em Reproduction of the linear functions}:
\begin{align}
\label{Eq:LinearReproduction}
\forall \M f\in\R^d,\ \forall \V x\in\R^d: \qquad \sum_{i=1}^N\langle \M f,\M r_i\rangle\,\varphi_i(\V x)=\langle \M f,\V x\rangle.
\end{align}
\item {\em Riesz stability}: for every $s=\sum_{i=1}^N c_i\varphi_i\in\CPWL_{\Spc F}$,
\begin{align}
\label{Eq:Riesz}
\max_{1\le i\le N}|c_i| \ \le\ \|s\|_{L_\infty(\mathbb S^{d-1})}\ \le\ \varsigma_{\Spc F}\,\max_{1\le i\le N}|c_i|,
\end{align}
where $\varsigma_{\Spc F}=\max_{\alpha}\sup\big\{\sum_{j=1}^d\tilde{\M r}_{\alpha,j}^{\Tr}\V \xi:\ \V \xi\in C_\alpha\cap\mathbb S^{d-1}\big\}$,
so that $\V c\mapsto s$ is a linear isomorphism from $(\R^N,\|\cdot\|_\infty)$ onto $\big(\CPWL_{\Spc F},\|\cdot\|_{L_\infty(\mathbb S^{d-1})}\big)$ with a condition number that does not exceed $\varsigma_{\Spc F}$.
\end{enumerate}
\end{proposition}

\begin{proof}
The key observation is that the dual (biorthogonal) system for the cell $C_\alpha$ is explicitly given by
\begin{align*}
\begin{bmatrix} \tilde{\M r}_{\alpha,1} & \cdots & \tilde{\M r}_{\alpha,d} \end{bmatrix}^{\Tr} = \M R_\alpha^{-1},
\end{align*}
the invertibility of $\M R_\alpha$ being guaranteed by the linear independence of the generating rays of the simplicial cell.

\medskip
\emph{Piecewise linearity and continuity}:
On any cell $C_\alpha$, the B-spline specified by \eqref{Eq:ConicBspline} is linear. To establish continuity across cell boundaries, let $C_\alpha$ and $C_\beta$ be two adjacent cells that share a common face. Let us suppose that $\M r_i$ is a vertex of this shared face, so that $\M r_i = \M r_{\alpha,k} = \M r_{\beta,m}$. For any point $\V x$ that lies on this face, the evaluations $\tilde{\M r}_{\alpha,k}^{\Tr}\V x$ and $\tilde{\M r}_{\beta,m}^{\Tr}\V x$ both extract the coordinate of $\V x$ along the $\M r_i$-direction and therefore coincide. Conversely, if $\M r_i$ does not belong to the shared face, then $\varphi_i(\V x) = 0$ on both sides of the boundary. Because the linear pieces match on all shared faces, global continuity is guaranteed, which confirms that $\varphi_i \in \CPWL_{\Spc F}$. The non-negativity and the homogeneity \eqref{Eq:BsplinePositivity}--\eqref{Eq:BsplineHomogeneity} are read off \eqref{Eq:ConicBspline}: on $C_\alpha$, the coordinate $\tilde{\M r}_{\alpha,k}^{\Tr}\V x$ of $\V x$ along $\M r_{\alpha,k}$ is non-negative and linear.

\medskip
\emph{Nodal property}:
Let $\M r_j$ be a ray. If $j=i$, then the choice of a cell $C_\alpha$ that contains $\M r_i=\M r_{\alpha,k}$ gives $\varphi_i(\M r_i) = \tilde{\M r}_{\alpha,k}^{\Tr}\M r_{\alpha,k} = 1$. If $j \neq i$ and the two rays share a cell $C_\alpha$ with $\M r_i=\M r_{\alpha,k}$ and $\M r_j=\M r_{\alpha,m}$, then $\varphi_i(\M r_j) = \tilde{\M r}_{\alpha,k}^{\Tr}\M r_{\alpha,m} = 0$. If they do not share a cell, then $\varphi_i(\M r_j)=0$ by definition. Hence, $\varphi_i(\M r_j) = \delta_{ij}$. The identification of $\operatorname{supp}(\varphi_i)$ as the union of the cells that contain $\M r_i$ follows directly from \eqref{Eq:ConicBspline}.

\medskip
\emph{Representation of conic splines and reproduction of the linear functions}:
Since the rays $\{\M r_{\alpha,j}\}_{j=1}^d$ form a basis of $\R^d$ with dual basis $\{\tilde{\M r}_{\alpha,j}\}_{j=1}^d$, any $\V x \in C_\alpha$ expands uniquely as $\V x = \sum_{j=1}^d c_j \M r_{\alpha,j}$ with $c_j=\tilde{\M r}_{\alpha,j}^{\Tr}\V x$ which, by \eqref{Eq:ConicBspline}, is the value at $\V x$ of the conic B-spline attached to the ray $\M r_{\alpha,j}$. Since $s \in \CPWL_{\Spc F}$ is linear on $C_\alpha$, we get
\begin{align}
\label{Eq:BsplineexpansionMin}
s(\V x) = s\Big(\sum_{j=1}^d c_j \M r_{\alpha,j}\Big)= \sum_{j=1}^d c_j\, s(\M r_{\alpha,j})= \sum_{i:\,\M r_i\in C_\alpha} s(\M r_i)\,\varphi_i(\V x).
\end{align}
The $(N-d)$ remaining B-splines, whose rays do not belong to $C_\alpha$, vanish on $C_\alpha$ by \eqref{Eq:ConicBspline}, so that the sum in \eqref{Eq:BsplineexpansionMin} extends to all $N$ rays, which is \eqref{Eq:NodalExpansion} on $C_\alpha$. As the cells of the complete fan $\Spc F$ cover $\R^d$, \eqref{Eq:NodalExpansion} holds for all $\V x\in\R^d$, and it shows that $s$ is uniquely determined by its nodal values.
Item 3 is \eqref{Eq:NodalExpansion} applied to the linear function $\V x\mapsto\langle \M f,\V x\rangle$, which is a member of $\CPWL_{\Spc F}$ with nodal values $\langle \M f,\M r_i\rangle$.

\medskip
\emph{Riesz stability}:
The lower bound in \eqref{Eq:Riesz} is the nodal property, which gives $c_i=s(\M r_i)$ with $\M r_i\in\mathbb S^{d-1}$. For the upper bound, let $\V \xi\in C_\alpha\cap\mathbb S^{d-1}$. Only the $d$ B-splines that are attached to the rays of $C_\alpha$ are non-zero there, so that \eqref{Eq:ConicBspline} and \eqref{Eq:BsplinePositivity} give $|s(\V \xi)|\le\max_i|c_i|\sum_{j=1}^d\tilde{\M r}_{\alpha,j}^{\Tr}\V \xi\le \varsigma_{\Spc F}\max_i|c_i|$.
\end{proof}

Since each dual vector $\tilde{\M r}_{\alpha,j}$ is orthogonal to the facet of $C_\alpha$ that is spanned by all its generating rays except $\M r_{\alpha,j}$, it is identified as the inward-pointing normal of that facet, with a scaling that is imposed by the biorthogonality condition \eqref{Eq:Biorthogonality}. This enables the geometric construction of the conic B-splines by continuously gluing together the local linear forms $\tilde{\M r}_{\alpha,j}^{\Tr}\V x$ across the adjacent simplicial cones of the fan. 
The reproduction of the linear functions, Item 3, is what gives the conic spline interpolant its second-order accuracy in Section~\ref{Sec:Regularity}.

Proposition~\ref{Prop:ConicBsplines} is the conic counterpart of the characterization of the CPWL functions on a simplicial triangulation of $\R^d$ by the hat functions of that triangulation \cite{Goujon2023stable}, with the Riesz constants depending on the shape of the simplices alone. The conic B-splines are the hat functions of the spherical complex of Section~\ref{Sec:RegularityIndices}, extended homogeneously to $\R^d$. Since the conic B-splines, unlike the hat functions, do not sum to one, the price of that extension is 
the factor $\varsigma_{\Spc F}$ in \eqref{Eq:Riesz}.

\subsection{Representation of Polyhedral Quasi-seminorms}
\label{Sec:PolyRepresentation}

\begin{theorem}[Characterization of the polyhedral CH functionals]
\label{Theo:PolyQuasiSemiNorms}
Let $q:\R^d\to[0,\infty)$ be an even functional. Then, the following statements are equivalent.
\begin{enumerate}
\item $q$ is a polyhedral CH functional; that is, $q$ is continuous piecewise-linear (CPWL) and absolutely homogeneous.
\item There exist a complete polyhedral fan $\Spc F$ in $\R^d$ and vectors $\{\M f_\alpha\}_{C_\alpha\in\Spc F}\subset\R^d$ such that
\begin{align*}
q(\V x)=\langle \M f_\alpha,\V x\rangle, \qquad \V x\in C_\alpha,
\end{align*}
with the compatibility condition on shared faces
\begin{align*}
\langle \M f_\alpha,\V x\rangle = \langle \M f_\beta,\V x\rangle, \qquad \V x\in C_\alpha\cap C_\beta.
\end{align*}
\item There exist a complete symmetric simplicial fan $\widetilde{\Spc F}$ with normalized rays $\Spc R=\{\M r_1,\dots,\M r_N\}\subset\mathbb S^{d-1}$ and associated conic B-splines $\{\varphi_n\}_{n=1}^N\subset \CPWL_{\widetilde{\Spc F}}$ such that
\begin{align}
\label{Eq:qBsplineExpansion}
q(\V x) = \sum_{n=1}^N q(\M r_n)\,\varphi_n(\V x).
\end{align}
\item The unit ball $B_q=\{\V x\in\R^d: q(\V x)\le 1\}$ is a symmetric star-shaped polyhedron that contains the origin in its interior.
\end{enumerate}
Such a functional is a {\em polyhedral quasi-seminorm} if and only if its null space $\Spc N_q$ is a linear subspace of $\R^d$; see Proposition~\ref{Prop:ConicSplineAdmissible}. It is a {\em polyhedral quasi-norm} if and only if the nodal values $q(\M r_n)$ of Item 3 are all strictly positive, in which case $B_q$ is bounded and hence is a symmetric star-shaped polytope whose boundary vertices are given by
\begin{align}
\label{Eq:Boundaryvertices}
\M v_n=\frac{\M r_n}{q(\M r_n)}.
\end{align}
The symmetry of $q$ induces the following structural relations:
\begin{itemize}
\item for every ray $\M r_n$, there exists a ray $\M r_{n'}=-\M r_n$ such that
\begin{align*}
q(\M r_{n'})=q(\M r_n), \qquad \varphi_{n'}(\V x)=\varphi_n(-\V x);
\end{align*}
\item for every cone $C_\alpha$ of $\Spc F$ or $\widetilde{\Spc F}$, the opposite cone $-C_\alpha=C_{-\alpha}$ also belongs to the fan, with the corresponding coefficient vector $\M f_{-\alpha}=-\M f_\alpha$.
\end{itemize}
\end{theorem}
\begin{proof}
We prove the cycle of implications $(1)\implies(2)\implies(3)\implies(4)\implies(1)$, and then establish the additional properties that are associated with positive-definiteness.

\medskip
\emph{$(1)\Rightarrow(2)$}:
Let us assume that $q$ is CPWL and absolutely homogeneous. The CPWL property implies the existence of a finite collection of closed convex polyhedral domains $C_\alpha$---the maximal regions of linearity of $q$---whose union is $\R^d$ and on each of which $q$ is affine; see Proposition~\ref{Prop:PartitionCPWL} in Appendix~\ref{App:Polyhedral}. Because $q(\lambda \V x) = \lambda q(\V x)$ for all $\lambda \ge 0$, the function $q$ is positively homogeneous. Consequently, if a vector $\V x$ belongs to a specific linearity domain, then the entire ray $\{\lambda \V x : \lambda \ge 0\}$ resides within that same domain. Moreover, since $q(\V 0)=0$, the constant intercept of each affine piece necessarily vanishes. Thus, $q$ is purely linear on each domain, which takes the form $q(\V x) = \langle \M f_\alpha, \V x \rangle$.

Due to the ray-inclusion property above, each linearity domain $C_\alpha$ is a polyhedral cone. Since the collection of these cones covers $\R^d$ and any two distinct cones intersect at most along a common face, they constitute a complete polyhedral fan $\Spc F$. The continuity of $q$ across the entire space forces the linear pieces to match on shared boundaries, so that $\V x \in C_\alpha \cap C_\beta \implies \langle \M f_\alpha, \V x \rangle = \langle \M f_\beta, \V x \rangle$. Finally, since $q$ is even, the fan is centrally symmetric: if $\V x \in C_\alpha$, then $-\V x \in -C_\alpha$ and $q(\V x)=q(-\V x)$. Therefore, $\langle \M f_\alpha,\V x\rangle = \langle \M f_{-\alpha},-\V x\rangle$, which implies that $\M f_{-\alpha}=-\M f_\alpha$.

\medskip
\emph{$(2)\Rightarrow(3)$}:
By invoking Proposition~\ref{Prop:FanSimplicial}, we refine the symmetric polyhedral fan $\Spc F$ into a complete symmetric simplicial fan $\widetilde{\Spc F}$ whose normalized rays are $\Spc R = \{\M r_1, \dots, \M r_N\} \subset \mathbb{S}^{d-1}$. When a cone of $\Spc F$ is not pointed---which happens, for instance, for the two half-planes of the linearity fan of $q(\V x)=|x_2|$ in $\R^2$---the refinement introduces new rays, which is harmless here, since a refinement of $\Spc F$ is all that we need. The fact that $\widetilde{\Spc F}$ is a refinement of $\Spc F$ implies that $q \in \CPWL_{\widetilde{\Spc F}}$, which enables us to apply Proposition~\ref{Prop:ConicBsplines} to obtain the conic B-spline expansion \eqref{Eq:qBsplineExpansion}. The symmetry of the fan, where $\M r_n \in \Spc R \implies -\M r_n \in \Spc R$, coupled with the evenness of $q$, translates into the stated structural conditions, which are characteristic of a symmetric B-spline expansion.

\medskip
\emph{$(3)\Rightarrow(4)$}:
Any function $q$ that is defined by the symmetric B-spline expansion \eqref{Eq:qBsplineExpansion} is CPWL and symmetric by construction. Due to the non-negativity of the expansion coefficients $q(\M r_n)$, it also inherits the homogeneity and the non-negativity of the constituent conic B-splines, see \eqref{Eq:BsplinePositivity}--\eqref{Eq:BsplineHomogeneity}. Hence, $q \in \CH(\R^d\to\R_+)$ and, by the same argument as in the proof of Theorem~\ref{Theo:SemiQuasiNorm}, its unit ball $B_q$ is symmetric, closed, star-shaped, and contains $\V 0$ as an interior point.

To show that $B_q$ is a polyhedron, we consider an individual simplicial cone $C_\alpha$ of the fan $\widetilde{\Spc F}$ that is associated with \eqref{Eq:qBsplineExpansion}. Within this cone, $q(\V x) = \langle \M f_\alpha, \V x \rangle$ is a single linear function, so that the intersection of $B_q$ with the cone is $\{\V x \in C_\alpha : \langle \M f_\alpha, \V x \rangle \le 1\}$, which is the solution set of finitely many linear inequalities and hence a polyhedron. Since $B_q$ is the finite union of these polyhedral sectors over all the cones of $\widetilde{\Spc F}$ and their boundaries align continuously, $B_q$ is a polyhedron.

\medskip
\emph{$(4)\Rightarrow(1)$}:
Let us assume that (4) holds. Because $B_q$ is a symmetric star-shaped polyhedron, its boundary $\partial B_q$ is composed of a finite number of flat polyhedral facets. We define a polyhedral fan $\Spc F$ by constructing the cones from the origin to each facet of $B_q$. The gauge representation $q(\V x) = \inf \{ t > 0 : \V x \in t B_q \}$ means that, for any vector $\V x$, the value $q(\V x)$ is determined by the distance from the origin to the boundary $\partial B_q$ along the ray that contains $\V x$. Within a specific cone $C_\alpha$ that faces a single flat facet, the boundary facet lies on a hyperplane $\langle \M f_\alpha, \V x \rangle = 1$. Therefore, for any $\V x \in C_\alpha$, the scaling factor that is required to push $\V x$ to the boundary is exactly given by the linear evaluation $\langle \M f_\alpha, \V x \rangle$, so that $q(\V x) = \langle \M f_\alpha, \V x \rangle$ inside $C_\alpha$. Because $B_q$ is a finite polyhedron, there are finitely many such cones, which proves that $q$ is CPWL. Absolute homogeneity follows from the symmetry $B_q = -B_q$ of the polyhedral boundary.

\medskip
\emph{Positive-definiteness}:
The positive-definiteness of $q$ is equivalent to the strict positivity of all the nodal values $q(\M r_n)$, since $q$ vanishes on the ray $\M r_n$ if $q(\M r_n)=0$ and since \eqref{Eq:qBsplineExpansion} is a non-negative combination of the conic B-splines, which are positive in the interior of their support. It is also equivalent to the boundedness of $B_q$, while a bounded polyhedron is by definition a {\em polytope}. Under this configuration, the normalized rays $\M r_n$ point toward the boundary of the polytope. The boundary vertices along these ray directions are found by scaling each generator $\M r_n$ so that its functional value equals one. By setting $\M v_n = \lambda \M r_n$, we solve
\begin{align*}
q(\M v_n) = 1 \implies \lambda \, q(\M r_n) = 1 \implies \lambda = \frac{1}{q(\M r_n)},
\end{align*}
which yields \eqref{Eq:Boundaryvertices}.
\end{proof}

It remains to decide which of these functionals are quasi-seminorms. In the polyhedral setting, the answer is combinatorial and it is read off the zero pattern of the nodal values.

\begin{proposition}[Admissibility of a conic spline]
\label{Prop:ConicSplineAdmissible}
Let $\widetilde{\Spc F}$ be a complete symmetric simplicial fan with rays $\{\M r_n\}_{n=1}^N$, let $c_n\ge0$ be nodal values with $c_{n'}=c_n$ whenever $\M r_{n'}=-\M r_n$, and let $q=\sum_{n=1}^N c_n\varphi_n$ be the associated conic spline. Then, $q$ is a polyhedral quasi-seminorm if and only if its null space $\Spc N_q=\{\V x \in \R^d: q(\V x)=0\}$ is a linear subspace of $\R^d$, which happens if and only if the rays with $c_n=0$ are the rays of a subfan of $\widetilde{\Spc F}$ whose support is a linear subspace $\Spc N$. In that case, $\Spc N_q=\Spc N$ and there exist $0<m_0\le M_0<\infty$ such that
\begin{align}
\label{Eq:CPWLDistanceComparison}
\forall \V x \in \R^d:\quad m_0 \dist(\V x,\Spc N) \le q(\V x)\le M_0 \dist(\V x,\Spc N).
\end{align}
In particular, the conic spline is a polyhedral quasi-norm whenever $c_n>0$ for every $n$.
\end{proposition}

\begin{proof}
The necessity of the condition is Theorem~\ref{Theo:SemiQuasiNorm}, Item 3. For the sufficiency, let us assume that $\Spc N_q$ is a linear subspace $\Spc N$. Since $q$ is CPWL and vanishes on $\Spc N$, it is Lipschitz-continuous with some constant $L$, which gives $q(\V x)=|q(\V x)-q(\V n)|\le L\|\V x-\V n\|$ for every $\V n \in \Spc N$ and hence the right-hand side of \eqref{Eq:CPWLDistanceComparison} with $M_0=L$. For the left-hand side, let us argue on a single cell $C_\alpha$ of the fan, on which $q(\V x)=\langle\M f_\alpha,\V x\rangle$, and let us consider the quotient $q(\V x)/\dist(\V x,\Spc N)$ on the compact set $K_\alpha=\{\V x \in C_\alpha: \|\V x\|_2=1,\ \V x\notin\Spc N\}$. If the infimum of this quotient over $K_\alpha$ were zero, then a minimizing sequence would admit a limit point $\V x_\star \in C_\alpha\cap\mathbb S^{d-1}$ with either $q(\V x_\star)=0$ and $\dist(\V x_\star,\Spc N)>0$, which contradicts $\Spc N_q=\Spc N$, or simply $\V x_\star \in \Spc N$. In this latter scenario, the linearity of $q$ on $C_\alpha$ and the linearity of $\V x\mapsto\dist(\V x,\Spc N)$ on the cone that is generated by $C_\alpha\cap\Spc N^\perp$ make the quotient constant along the rays that approach $\V x_\star$, so that its infimum is attained at a point where it is positive. The minimum over the finitely many cells of the fan gives $m_0>0$, and \eqref{Eq:CPWLDistanceComparison} then implies the quasi-triangle inequality by Item 3 of Theorem~\ref{Theo:SemiQuasiNorm}.
The last statement is the case $\{n: c_n=0\}=\varnothing$, in which the empty subfan forces $\Spc N_q=\{\V 0\}$.
\end{proof}

The subfan condition is void as long as all the nodal values are positive, which is the configuration of the quasi-norms, and it is easy to monitor otherwise: what it forbids is an isolated zero, such as the pair of vanishing diagonal values of Example~\ref{Ex:NotQuasiSeminorm}.
That said, a design that requires a null space is better served by Theorem~\ref{Theo:QuasiSeminormStructure}, Item 1(iii), where the null space is carried by a linear map and the conic spline stays positive.

Among the four representations in Theorem~\ref{Theo:PolyQuasiSemiNorms}, only Items 3 and 4 are genuinely constructive; the others primarily establish analytical and geometric properties.

The B-spline representation (Item 3) is the most direct description. It implicitly guarantees global continuity by hardcoding the underlying domain partition into a complete simplicial fan $\widetilde{\Spc F}$. This structural constraint allows for the seamless parametrization of the whole of $\CPWL_{\widetilde{\Spc F}}\cap\CH(\R^d\to\R_+)$ simply by prescribing the ray values $q_n = q(\M r_n) \ge 0$ in \eqref{Eq:qBsplineExpansion}, subject only to the point-symmetry requirement $q_n = q_{n'}$ whenever $\M r_n = -\M r_{n'}$. On the open cone $\{q_n>0\}$, every such choice is a polyhedral quasi-norm, so that the parametrization is free of constraints where it matters; on the boundary of that cone, the admissibility is decided by Proposition~\ref{Prop:ConicSplineAdmissible}.

Item 4 offers an indirect characterization of the quasi-(semi)norm via the geometry of its unit ball $B_q$. The construction of $q$ then reduces to the design of an admissible symmetric polyhedral body whose gauge functional $q = \Gamma_{B_q}$ recovers the quasi-seminorm.

When $q$ specializes to a valid quasi-norm, these two constructive forms are coupled through the projective relation \eqref{Eq:Boundaryvertices}. This formula acts as a bridge: it dictates how to compute the vertices of the bounded polytope $B_q$ from the ray coefficients of Item 3.

By comparison, the classical description in Item 1 is of limited value for practical design. While it asserts the existence of linear pieces, it provides no structural mechanism to enforce the boundary compatibility and the continuity conditions across the shared faces of adjacent polyhedral domains.

\section{Regularity and Approximation by Conic Splines}
\label{Sec:Regularity}

\subsection{Regularity Indices of a Quasi-seminorm}

\label{Sec:RegularityIndices}
Let us assume that $q:\R^d\to[0,+\infty)$ is a continuous quasi-(semi)norm on $\R^d$, so that it is absolutely homogeneous. We propose to quantify its first- and second-order regularity by the computation of its gradient $\nabla\{q\}$ and of its Hessian $\Op H\{q\}$ over $\R^d$. Since $q$ is homogeneous, it is sufficient to assess this regularity on the unit sphere $\mathbb{S}^{d-1}$.

Inspired by the concepts of first- and second-order total variation \cite{Rudin1992,Ambrosio2000,Aziznejad2023a}, we introduce the corresponding spherical analogues, which are obtained by integrating a suitable differential norm over $\mathbb{S}^{d-1}$. Specifically, under the assumption that $q$ is differentiable on $\R^d\backslash\{\V 0\}$, we define
\begin{align}
\label{Eq:Lipq}
\Lip_{\mathbb{S}^{d-1}}(q)&= \sup_{\V \xi \in \mathbb{S}^{d-1}} \|\nabla\{q\}(\V \xi)\|_{2},\\
\label{Eq:TVq}
\TV_{\mathbb{S}^{d-1}}(q)&= \int_{\mathbb{S}^{d-1}} \|\nabla\{q\}(\V \xi)\|_{2} \; \dint \Spc H^{d-1}(\V \xi),
\end{align}
where $\Spc H^{d-1}$ stands for the $(d-1)$-dimensional Hausdorff (or surface) measure on $\mathbb S^{d-1}$. The symbol $\Lip$ in \eqref{Eq:Lipq} stands for Lipschitz, as the so-defined quantity yields the Lipschitz constant of $q$ over $\R^d \backslash\{\V 0\}$, which is the $L_\infty$-norm of the gradient. Indeed, since $q$ is $1$-homogeneous, its gradient is $0$-homogeneous in the sense that $\nabla\{q\}(\lambda\V x)=\lambda^0\,\nabla\{q\}(\V x)=\nabla\{q\}(\V x)$ for $\lambda>0$, which makes the supremum over the compact unit sphere equivalent to the supremum over all of $\R^d\backslash\{\V 0\}$.

The second-order indices involve the Hessian $\Op H\{q\}$ instead of the gradient. When $q$ is twice differentiable away from the origin, they are
\begin{align}
\label{Eq:HLip}
\HLip_{\mathbb{S}^{d-1}}(q)&= \sup_{\V \xi \in \mathbb{S}^{d-1}} \|\Op H\{q\}(\V \xi)\|_{S_\infty},\\
\label{Eq:HTVq}
\HTV_{\mathbb{S}^{d-1}}(q)&= \int_{\mathbb{S}^{d-1}} \|\Op H\{q\}(\V \xi)\|_{S_1}\, \dint \Spc H^{d-1}(\V \xi),
\end{align}
where $\|\M W\|_{S_p}=\big\|\big(\sigma_i(\M W)\big)\big\|_{\ell_p}$ denotes the Schatten-$p$ norm of the symmetric matrix $\M W \in \R^{d \times d}$ with singular values $\sigma_1(\M W)\ge\dots \ge \sigma_d(\M W)\ge 0$; in particular, $\|\M W\|_{S_\infty}=\sigma_1(\M W)$ and $\|\M W\|_{S_1}=\sum_i\sigma_i(\M W)$. Since the sphere is compact, $\HLip_{\mathbb{S}^{d-1}}(q)<\infty$ implies $\HTV_{\mathbb{S}^{d-1}}(q)<\infty$, while the reverse implication does not hold.

The index \eqref{Eq:HTVq} is the one that stays finite on the polyhedral members of the class, whose Hessian is a measure rather than a function. The definition that covers them is obtained by duality, as in \cite{Aziznejad2023a}.

\begin{definition}[Spherical Hessian variation]
\label{Def:HTV}
Let $q:\R^d\to\R_+$ be continuous and $1$-homogeneous, with $d\ge2$. Its {\em spherical Hessian variation} is
\begin{align}
\label{Eq:HTVgeneral}
\HTV_{\mathbb{S}^{d-1}}(q)=(d-1)\,\sup_{\V \Phi}\int_{B_{\ell_2}}q(\V x)\,\sum_{i,j=1}^d\partial_i\partial_j\Phi_{ij}(\V x)\, \dint\V x,
\end{align}
where $B_{\ell_2}$ is the Euclidean unit ball and where the supremum is taken over the symmetric matrix-valued test fields $\V \Phi\in C^\infty_c(B_{\ell_2};\R^{d\times d})$, compactly supported in the interior of $B_{\ell_2}$, with $\|\V \Phi(\V x)\|_{S_\infty}\le1$ for every $\V x$. Equivalently, $\HTV_{\mathbb{S}^{d-1}}(q)=(d-1)\,|\Op H\{q\}|(B_{\ell_2})$ whenever the distributional Hessian $\Op H\{q\}$ is a matrix-valued Radon measure, where $|\Op H\{q\}|$ denotes its total-variation measure with respect to the Schatten-$1$ norm.
\end{definition}

Definition~\ref{Def:HTV} is compatible with \eqref{Eq:HTVq}. Indeed, if $q$ is twice differentiable away from the origin, then its Hessian is $(-1)$-homogeneous, so that
\begin{align*}
(d-1)\int_{B_{\ell_2}}\|\Op H\{q\}\|_{S_1}\dint\V x&=(d-1)\int_0^1 r^{d-2}\dint r\int_{\mathbb S^{d-1}}\|\Op H\{q\}(\V \xi)\|_{S_1}\dint\Spc H^{d-1}(\V \xi)\\
&=\int_{\mathbb S^{d-1}}\|\Op H\{q\}(\V \xi)\|_{S_1}\dint\Spc H^{d-1}(\V \xi),
\end{align*}
which is \eqref{Eq:HTVq}. The same homogeneity shows that the value depends on the restriction of $q$ to the sphere alone.

\medskip
\emph{Evaluation for a polyhedral $q$}:
Let $\Spc F$ be a complete simplicial fan and let $q_{\Spc F}$ be absolutely homogeneous and CPWL with respect to $\Spc F$, so that $q_{\Spc F}(\V x)=\langle \M f_\alpha,\V x\rangle$ for $\V x\in C_\alpha$. The evaluation of the indices rests on the partition of the sphere that the fan induces.

\begin{definition}[Spherical complex associated with a fan]
\label{Def:SphericalComplex}
Let $\Spc F$ be a complete simplicial fan in $\R^d$. The \emph{spherical complex} associated with $\Spc F$ is
\begin{align*}
\Spc S = \bigl\{ F\cap\mathbb S^{d-1} :\; F \text{ is a nonzero face of a cone } C\in\Spc F \bigr\}.
\end{align*}
The elements of $\Spc S$ are called \emph{spherical faces}.
\end{definition}

The {\em spherical cells} of $\Spc S$ are the $(d-1)$-dimensional spherical simplices $\Sigma_\alpha=C_\alpha\cap\mathbb S^{d-1}$, where $C_\alpha$ ranges over the cells of $\Spc F$. Two adjacent cells share the {\em spherical facet} $F_{\alpha\beta}=\Sigma_\alpha\cap\Sigma_\beta$, which is a $(d-2)$-dimensional spherical simplex, and $\adj(\alpha)$ denotes the set of the indices of the cells that are adjacent to $\Sigma_\alpha$.

The gradient of $q_{\Spc F}$ is the constant $\M f_\alpha$ on the interior of each cone, hence on each spherical cell, so that \eqref{Eq:TVq} reduces to
\begin{align}
\label{Eq:TVCPWL}
\TV_{\mathbb{S}^{d-1}}(q_{\Spc F})=\sum_{\alpha \in A} \|\M f_\alpha\|_2\, \Spc H^{d-1}(\Sigma_\alpha),
\end{align}
where $\Spc H^{d-1}(\Sigma_{\alpha})$ is the spherical surface of $\Sigma_{\alpha}$. The Hessian of $q_{\Spc F}$ vanishes inside every cone and is carried by the facets. Across the facet $F_{\alpha\beta}$, the continuity of $q_{\Spc F}$ forces the gradient jump $\M f_\alpha-\M f_\beta$ to be parallel to the unit normal $\V n_{\alpha\beta}$ of the facet, so that $\Op H\{q_{\Spc F}\}$ is $\|\M f_\alpha-\M f_\beta\|_2\,\V n_{\alpha\beta}\V n_{\alpha\beta}^\Tr$ times the surface measure of the cone over $F_{\alpha\beta}$. The rank-one factor has unit Schatten-$1$ norm, and the cone over $F_{\alpha\beta}$ meets $B_{\ell_2}$ in a set of surface measure $\Spc H^{d-2}(F_{\alpha\beta})/(d-1)$. Definition~\ref{Def:HTV} therefore gives
\begin{align}
\label{Eq:HTVCPWL}
\HTV_{\mathbb{S}^{d-1}}(q_{\Spc F})=\frac{1}{2}\sum_{\alpha \in A} \sum_{\beta \in \adj(\alpha)}
\|\M f_\alpha-\M f_\beta\|_2\, \Spc H^{d-2}(F_{\alpha\beta}),
\end{align}
where the factor $\frac12$ corrects the double counting of each facet in the double sum. This is the exact spherical analogue of the formula for $\HTV(f)$ of a CPWL function $f$ on a polyhedral partition of $\R^d$ \cite[Theorem 4.1]{Aziznejad2023a}, the dimension of the facets being one less here.

\begin{figure}[t]
\centering
\includegraphics[width=\textwidth]{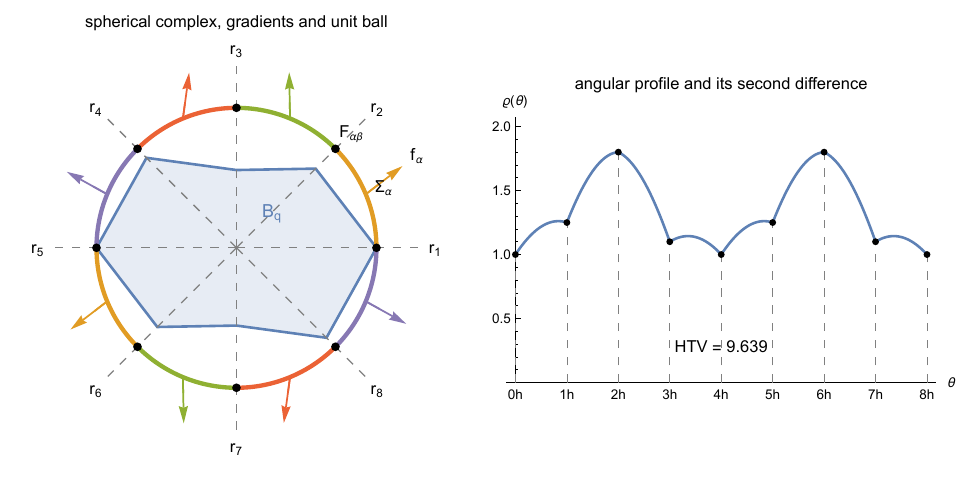}
\caption{\small The spherical complex and the evaluation of \eqref{Eq:HTVCPWL} for $d=2$. Left: a symmetric fan with eight rays $\M r_1,\dots,\M r_8$ (dashed), the spherical cells $\Sigma_\alpha$ (arcs of the unit circle), the breakpoints $F_{\alpha\beta}=\Sigma_\alpha\cap\Sigma_\beta$ (dots), the constant gradient $\M f_\alpha$ of $q$ on each cone, and the star-shaped unit ball $B_q$ of the polyhedral quasi-norm $q$ with nodal values $q(\M r_i)\in\{1,1.25,1.8,1.1\}$, repeated antipodally. The ball is not convex: at the two vertices of smallest radius, which lie on the rays with the nodal value $1.8$, the interior angle exceeds $\pi$ and the boundary is indented. Right: the angular profile $\varrho(\theta)=q(\V\xi(\theta))=\langle\M f_\alpha,\V\xi(\theta)\rangle$, which is a sinusoidal arc on each cell; the Hessian measure of $q$ is concentrated on the breakpoints, with the weights $\|\M f_\alpha-\M f_\beta\|_2$, and $\HTV_{\mathbb S^1}(q)$ is their sum.}
\label{Fig:SphericalComplex2D}
\end{figure}

Figure~\ref{Fig:SphericalComplex2D} illustrates the construction for $d=2$, where the spherical cells are the arcs between consecutive rays, the facets $F_{\alpha\beta}$ are the breakpoints, and $\Spc H^{0}$ is the counting measure, so that \eqref{Eq:HTVCPWL} reduces to the sum of the jumps $\|\M f_\alpha-\M f_\beta\|_2$ of the gradient over the breakpoints.

\subsection{Approximation by Conic Splines}

\label{Sec:Approximation}
We now investigate the possibility of approximating any quasi-seminorm $q$ by a CPWL one to an arbitrary degree of precision. To do so, we formulate the task as a conic spline interpolation problem and derive the corresponding error bounds under the assumption that $q$ is twice differentiable. The estimates involve two standard geometric requirements on the sequence of fans, which we state once and for all.

\begin{definition}[Shape-regular and quasi-uniform fans]
\label{Def:ShapeRegular}
Let $\Spc F_n$ be a complete simplicial fan of $\R^d$ with normalized rays $\M r_1,\dots,\M r_n \in \mathbb S^{d-1}$. Its spherical cells are $\Sigma_\alpha=C_\alpha\cap\mathbb S^{d-1}$ and the associated Euclidean simplices are $T_\alpha=\Conv\{\M r_{\alpha,1},\dots,\M r_{\alpha,d}\}$. The fan is said to be
\begin{itemize}
\item {\em shape-regular} with constant $\kappa\ge1$ if, for every cell, the ratio of the circumradius of $T_\alpha$ to its inradius, both taken relative to the affine hull of $T_\alpha$, does not exceed $\kappa$;
\item {\em quasi-uniform} with constant $\nu\ge1$ if the ratio of the largest to the smallest cell diameter does not exceed $\nu$.
\end{itemize}
We write $h_n=\max_\alpha \diam(\Sigma_\alpha)$ for the maximal spherical cell diameter. A family of fans is shape-regular (resp., quasi-uniform) if the constant may be chosen independently of $n$.
\end{definition}

\begin{theorem}[Approximation of $1$-homogeneous functions by conic splines]
\label{Theo:ConicSplineApprox}
Let $q:\R^d\to\R$ be $1$-homogeneous and $C^2$ on $\R^d\backslash\{\V{0}\}$. Let $\Spc{F}_n$ be a shape-regular simplicial fan with rays $\M{r}_1,\ldots,\M{r}_n\in\mathbb S^{d-1}$ and maximal spherical cell diameter $h_n$. We construct the spline interpolant
\begin{align}
\label{Eq:SplineInterpolant}
q_n(\V x) = \sum_{i=1}^n q(\M{r}_i)\varphi_i(\V x),
\end{align}
where the functions $\varphi_i$ are the conic B-splines induced by $\Spc{F}_n$. Then, the following holds.
\begin{enumerate}
\item {\em Uniform bound}:
\begin{align}
\label{Eq:ApproxBound}
\sup_{\V \xi \in \mathbb{S}^{d-1}} |q(\V \xi)-q_n(\V \xi)|\le
C_{\Spc{F}}\,h_n^2\, \HLip_{\mathbb S^{d-1}}(q),
\end{align}
where $\HLip_{\mathbb S^{d-1}}(q)$ is the second-order smoothness index specified by \eqref{Eq:HLip} and where the constant $C_{\Spc{F}}$ depends solely on the shape-regularity of the fan.
\item {\em Rate}: if the rays are quasi-uniform on $\mathbb S^{d-1}$, then $h_n\asymp n^{-1/(d-1)}$, which yields
\begin{align}
\label{Eq:ApproxRate}
\|q-q_n\|_{L_\infty(\mathbb{S}^{d-1})} \le C'_{\Spc{F}}\, \frac{1}{n^{2/(d-1)}}\, \HLip_{\mathbb S^{d-1}}(q).
\end{align}
\item {\em Sharpness}: the rate \eqref{Eq:ApproxRate} cannot be improved, and the obstruction is already present for the most regular member of the class. For the Euclidean norm $q=\|\cdot\|_2$ with $d\ge2$, there exists a constant $c_d>0$ such that, for every complete simplicial fan $\Spc F$ with at most $n$ rays and every conic spline $s\in\CPWL_{\Spc F}$,
\begin{align}
\label{Eq:LowerBound}
\|q-s\|_{L_\infty(\mathbb S^{d-1})}\ \ge\ c_d\, n^{-2/(d-1)}.
\end{align}
The approximation of a quasi-norm by conic splines therefore saturates at the order $h_n^2$, however smooth the quasi-norm is.
\end{enumerate}
\end{theorem}

\begin{proof}
\emph{Constants provided by shape-regularity}:
Let $T_\alpha=\Conv(\M r_{\alpha,1},\dots,\M r_{\alpha,d})$ be the Euclidean simplex of the cell $C_\alpha$. Shape-regularity provides constants that depend on $\kappa$ and on $d$ alone: $C_1$ with $\diam(T_\alpha)\le C_1h_n$; $c_2$ with $\Spc H^{d-1}(\Sigma_\alpha)\le c_2h_n^{d-1}$; and $c>0$ with $\min_{\V y\in T_\alpha}\|\V y\|_2\ge c$, the last one because the vertices of $T_\alpha$ lie on the unit sphere and a simplex whose circumradius-to-inradius ratio is bounded cannot have a face that passes arbitrarily close to the center of that sphere.
The normalization factor $\varsigma_{\Spc F}$ of Proposition~\ref{Prop:ConicBsplines}, Item 4, depends on the geometry of the fan alone, and it is bounded for the same reason: for $\V \xi\in C_\alpha\cap\mathbb S^{d-1}$, the sum of the barycentric coefficients $\sum_j\tilde{\M r}_{\alpha,j}^\Tr\V \xi$ is the factor $\varsigma$ with $\V \xi=\varsigma\V z$ and $\V z\in T_\alpha$, so that $\varsigma=1/\|\V z\|_2\le1/c$ and $\varsigma_{\Spc F}\le 1/c$. It is therefore bounded independently of $n$ for a shape-regular family, which is what makes it a legitimate constant in the estimates below. All the constants that are denoted by $C_{\Spc F}$, $C'_{\Spc F}$, $C_1$, $c$, and $c_2$ in this section are of this kind.

\medskip
\emph{Reduction to a flat simplex}:
Let $C=\Cone(\M{r}_{i_1},\ldots,\M{r}_{i_d})$ be a simplicial cone of the fan and let
\begin{align*}
T=\Conv(\M{r}_{i_1},\ldots,\M{r}_{i_d})
\end{align*}
be the associated $(d-1)$-simplex. For any $\V \xi\in C\cap\mathbb S^{d-1}$, there exist unique coefficients $\alpha_1,\ldots,\alpha_d\ge 0$ such that $\V \xi=\sum_{j=1}^{d}\alpha_j\M{r}_{i_j}$. Since $\|\V \xi\|_2=1$, not all coefficients vanish, and we set
\begin{align*}
\varsigma=\sum_{j=1}^{d}\alpha_j>0, \qquad \V z=\frac{\V \xi}{\varsigma} =\sum_{j=1}^{d}\frac{\alpha_j}{\varsigma}\,\M{r}_{i_j},
\end{align*}
with the guarantee that $\V z\in T$ because the coefficients $\alpha_j/\varsigma$ are non-negative and sum to unity, and with $\varsigma\le\varsigma_{\Spc F}$ by \eqref{Eq:Riesz}.

By construction, $q_n$ is the positively homogeneous extension of the affine interpolant $\Op I_T\{q\}$ on $T$ specified in Lemma~\ref{Lem:AffineInterp} below. Hence, we can utilize the homogeneity of both $q$ and $q_n$ to write
\begin{align*}
q(\V \xi)-q_n(\V \xi) = \varsigma\bigl(q(\V z)-\Op I_T\{q\}(\V z)\bigr).
\end{align*}
\medskip
\emph{Bounding the local interpolation error}:
We apply Lemma~\ref{Lem:AffineInterp}, which yields
\begin{align}
\label{Eq:ErrBound1}
|q(\V z)-\Op I_T\{q\}(\V z)| \le \frac12\,\diam(T)^2 \sup_{\V y\in T}\|\Op H \{q\}(\V y)\|_{S_p}.
\end{align}
We recall that $\diam(T)\le C_1 h_n$ and that $\|\V y\|_2\ge c$ for all $\V y\in T$. Since $q$ is $1$-homogeneous, its Hessian satisfies
\begin{align}
\label{Eq:HessianHomogeneity}
\Op H \{q\}(\lambda \V y) = \lambda^{-1}\Op H \{q\}(\V y), \qquad \lambda>0.
\end{align}
Consequently, we can write the bound as
\begin{align*}
\|\Op H \{q\}(\V y)\|_{S_p}
&=
\|\V y\|_2^{-1} \left\| \Op H \{q\}\!\left(\frac{\V y}{\|\V y\|_2}\right) \right\|_{S_p} \\
&\le
c^{-1} \sup_{\V \xi \in \mathbb{S}^{d-1}} \|\Op H \{q\}(\V \xi)\|_{S_p},
\end{align*}
which implies that
\begin{align}
\sup_{\V y\in T}\|\Op H \{q\}(\V y)\|_{S_p}\le c^{-1} \sup_{\V \xi \in \mathbb{S}^{d-1}}\|\Op H \{q\}(\V \xi)\|_{S_p}.
\end{align}
We substitute this expression back into the interpolation bound \eqref{Eq:ErrBound1} and take the supremum over all the cones of the fan. With the choice $p=\infty$, this operation yields the first estimate \eqref{Eq:ApproxBound} with $C_{\Spc{F}}= \frac{1}{2} c^{-1} C_1^2 \varsigma_{\Spc F}$.

\medskip
\emph{Asymptotic convergence rate}:
Under the assumption that the rays are quasi-uniform on $\mathbb S^{d-1}$, the classical covering estimate for the sphere gives $h_n\asymp n^{-1/(d-1)}$. We substitute this relation into the first inequality, which directly yields the announced rate \eqref{Eq:ApproxRate}.

\medskip
\emph{Sharpness}:
Let $q=\|\cdot\|_2$ and let $s\in\CPWL_{\Spc F}$ satisfy $\|q-s\|_{L_\infty(\mathbb S^{d-1})}\le\varepsilon<1$ on a fan with at most $n$ rays. Then, $s$ is non-negative and $1$-homogeneous, and its unit ball $B_s$ is a star-shaped polytope with at most $n$ vertices, namely the points $\M r_i/s(\M r_i)$. The two-sided bound $(1-\varepsilon)\le s(\V \xi)\le(1+\varepsilon)$ on the sphere gives $(1+\varepsilon)^{-1}B_{\ell_2}\subseteq B_s\subseteq(1-\varepsilon)^{-1}B_{\ell_2}$. Since $B_s$ need not be convex, we pass to the convex hull $P_s=\Conv\{\M r_i/s(\M r_i)\}$ of its vertices, which is a convex polytope with at most $n$ vertices. The latter satisfies $(1+\varepsilon)^{-1}B_{\ell_2}\subseteq B_s\subseteq P_s\subseteq(1-\varepsilon)^{-1}B_{\ell_2}$ 
because its vertices lie in the convex set $(1-\varepsilon)^{-1}B_{\ell_2}$. Hence, the Hausdorff distance between $P_s$ and the Euclidean ball $B_{\ell_2}$ is $O(\varepsilon)$. The classical lower bound for the approximation of a Euclidean ball by a convex polytope with $n$ vertices \cite{Gruber1993} then forces $\varepsilon\gtrsim n^{-2/(d-1)}$, which is \eqref{Eq:LowerBound}.
\end{proof}

\begin{lemma}[Control of the affine interpolation error on a simplex]
\label{Lem:AffineInterp}
Let $T=\Conv(\M{p}_0,\ldots,\M{p}_k)\subset\R^d$ be a $k$-simplex with $1\le k\le d$, and let $\lambda_0,\ldots,\lambda_k$ denote its barycentric coordinates.
Given $f\in C^2(V)$, where $V$ is an open neighborhood of $T$, we define the affine interpolant
\begin{align}
\Op I_T\{f\}(\V x) = \sum_{i=0}^{k} f(\M{p}_i)\lambda_i(\V x), \qquad \V x\in T.
\end{align}
Then, for every $p\in[1,\infty]$,
\begin{align}
\|f-\Op I_T\{f\}\|_{L_\infty(T)} \le \frac12\,\diam(T)^2 \sup_{\V x\in T}\|\Op H \{f\}(\V x)\|_{S_p}.
\end{align}
\end{lemma}

\begin{proof}
We start by writing Taylor's formula for $f(\M{p}_i)$ expanded at $\V x$ as
\begin{align*}
f(\M{p}_i) = f(\V x) +\nabla f(\V x)^{\Tr}(\M{p}_i-\V x) +\frac12 (\M{p}_i-\V x)^{\Tr} \Op H \{f\}(\V \xi_i) (\M{p}_i-\V x),
\end{align*}
for some $\V \xi_i\in[\V x,\M{p}_i]\subset T$. We then construct $\Op I_T\{f\}(\V x)=\sum_{i=0}^k\lambda_i f(\M{p}_i)$ by summing with the weights $\lambda_i$. The linear terms vanish because $\sum_{i=0}^k\lambda_i(\M{p}_i-\V x)=\V{0}$. Hence, we obtain
\begin{align}
\Op I_T\{f\}(\V x)-f(\V x) = \frac12 \sum_{i=0}^{k}\lambda_i (\M{p}_i-\V x)^{\Tr} \Op H \{f\}(\V \xi_i) (\M{p}_i-\V x).
\end{align}

\medskip
\emph{Matrix-norm optimization via duality}:
The key step is to rewrite the quadratic term as an inner product:
\begin{align*}
(\M{p}_i-\V x)^{\Tr} \Op H \{f\}(\V \xi_i) (\M{p}_i-\V x) = \left\langle \Op H \{f\}(\V \xi_i), \M{B}_i \right\rangle,
\end{align*}
where $\M{B}_i=(\M{p}_i-\V x)(\M{p}_i-\V x)^\Tr$. Since the matrix $\M{B}_i$ is symmetric and rank-one, its Schatten-$p'$ norm for any $p'\ge1$ is equal to its single non-zero eigenvalue, which implies that $\|\M{B}_i\|_{S_{p'}}=\|\M{p}_i-\V x\|_2^2$. By invoking H\"older's inequality for matrices (the duality relation for Schatten norms with $p,p' \in[1, \infty]$ and $\frac{1}{p} + \frac{1}{p'}=1$),
\begin{align}
|\langle \M{A},\M{B}\rangle| \le \|\M{A}\|_{S_p}\|\M{B}\|_{S_{p'}},
\end{align}
we then obtain the upper bound
\begin{align*}
\left| (\M{p}_i-\V x)^{\Tr} \Op H \{f\}(\V \xi_i) (\M{p}_i-\V x) \right| \le \|\Op H \{f\}(\V \xi_i)\|_{S_p} \|\M{p}_i-\V x\|_2^2.
\end{align*}
Since $\|\M{p}_i-\V x\|_2\le \diam(T)$ and $\sum_{i=0}^k \lambda_i = 1$ with $\lambda_i \ge 0$, it follows that
\begin{align}
|\Op I_T\{f\}(\V x)-f(\V x)| \le \frac12\diam(T)^2 \sup_{\V y\in T}\|\Op H \{f\}(\V y)\|_{S_p}.
\end{align}
We take the supremum over $\V x\in T$, which yields the desired estimate.
\end{proof}

Two consequences of the construction behind Theorem~\ref{Theo:ConicSplineApprox} are worth recording, because they are what turns an interpolation estimate into an approximation theorem. The first is that the interpolation is one-sided in the convex case: if $q$ is convex, then $q\le q_n$ pointwise, since $q_n$ is linear on each cell and agrees with $q$ at the rays, which reads $B_{q_n}\subseteq B_q$. The conic spline interpolant is then the gauge of the polytope inscribed in $B_q$ with vertices $\M r_i/q(\M r_i)$, which is the exact bridge to the polytopal approximation of convex bodies. The second one is that the interpolant is near-optimal: the combination of \eqref{Eq:ApproxRate} with \eqref{Eq:LowerBound} shows that, on a quasi-uniform fan, it achieves the best rate that any conic spline with $n$ rays can achieve, up to a constant that depends on $d$ and on the shape-regularity alone.

\paragraph{Relation to spherical splines.}
For $d=3$, the conic splines are the linear spherical splines of \cite{Alfeld1996a}, and Item 1 is consistent with the approximation theorem of Neamtu and Schumaker \cite[Theorem 1.1, Eq.~(5.6)]{Neamtu2004}, which bounds the error of the spherical splines of degree $m$ on a triangulation of $\mathbb S^2$ of mesh size $h$ by $C\,h^{m+1}\,|f|_{m+1,p,\mathbb S^2}$, with Sobolev $p$-seminorms that are defined through a homogeneous extension of $f$ and its projection onto a tangent plane. For $m=1$ and $p=\infty$, their seminorm is comparable to $\HLip_{\mathbb S^2}(q)$, and their proof reduces to the planar interpolation estimate in the same way as ours. 
What Theorem~\ref{Theo:ConicSplineApprox} contributes is the rate in an arbitrary dimension $d$, with the seminorm fixed by the $1$-homogeneity of $q$, together with the proof of its sharpness in Item 3, which has no counterpart in the theory of spherical splines. It is completed by Theorem~\ref{Theo:HTVApprox} below, which extends the second-order rate to measure-valued Hessians.

\subsection{Approximation under a Bounded Spherical Hessian Variation}
\label{Sec:ApproximationHTV}

Theorem~\ref{Theo:ConicSplineApprox} rests on the assumption that $q$ is twice differentiable away from the origin, which excludes the very objects that our parametrization produces: a polyhedral quasi-seminorm has a vanishing Hessian inside each cone and an infinite one across the facets, so that $\HLip_{\mathbb S^{d-1}}(q)=+\infty$ and the bound \eqref{Eq:ApproxBound} is vacuous. The index that remains finite is $\HTV_{\mathbb S^{d-1}}(q)$ of Definition~\ref{Def:HTV}, and the theorem of this section shows that it controls the $L_1$ error at the same second order, for every continuous quasi-seminorm with bounded spherical Hessian variation and without any differentiability assumption. Since $q$ is $1$-homogeneous, its Hessian measure is $(-1)$-homogeneous, which means that $|\Op H\{q\}|(\lambda B)=\lambda^{d-1}|\Op H\{q\}|(B)$ for every Borel set $B\subseteq\R^d\backslash\{\V 0\}$ and $\lambda>0$. In polar coordinates, $|\Op H\{q\}|$ therefore disintegrates as $r^{d-2}\dint r\otimes\mu_q$, where $\mu_q\le 0$ is a finite Borel measure on $\mathbb S^{d-1}$ with $\mu_q(\mathbb S^{d-1})=(d-1)|\Op H\{q\}|(B_{\ell_2})=\HTV_{\mathbb S^{d-1}}(q)$ by Definition~\ref{Def:HTV}; we call $\mu_q$ the {\em spherical Hessian measure} of $q$. For a Borel set $S\subseteq\mathbb S^{d-1}$ and $0<a<b$, the conic shell $K_S(a,b)=\{r\V \eta: \V \eta\in S,\ a<r<b\}$ then has the mass
\begin{align}
\label{Eq:ShellMass}
|\Op H\{q\}|\big(K_S(a,b)\big)=\mu_q(S)\int_a^b r^{d-2}\dint r .
\end{align}

\medskip
\emph{Preaveraged conic splines}:
Let $\chi:\R^d\to[0,\infty)$ be a fixed even $C^\infty$ window that is supported in the unit ball and normalized so that $\int_{\R^d}\chi=1$, and let $\chi_h=h^{-d}\chi(\cdot/h)$. For a continuous $1$-homogeneous $q$, we consider, next to the {\em sampled} conic spline $q_n$ of \eqref{Eq:SplineInterpolant}, the
{\em preaveraged} conic spline
\begin{align}
\label{Eq:SplineAveraged}
\tilde q_n(\V x)=\sum_{i=1}^n (q*\chi_{h_n})(\M r_i)\,\varphi_i(\V x),
\end{align}
where $(q*\chi_h)(\M r_i)=\int_{\R^d}\chi_h(\V y)\,q(\M r_i-\V y)\dint\V y$ is the local average of $q$ over the ball of radius $h$ centered at the ray $\M r_i$. Both splines are elements of $\CPWL_{\Spc F_n}$. Moreover, if $q$ is a quasi-seminorm and the fan is symmetric, then $q_n$ and $\tilde q_n$ are polyhedral CH functionals: their nodal values are non-negative, and they satisfy $q(-\M r_i)=q(\M r_i)$ and $(q*\chi_h)(-\M r_i)=(q*\chi_h)(\M r_i)$ because $q$ and $\chi$ are even, so that Theorem~\ref{Theo:PolyQuasiSemiNorms}, Item 3, applies.

The following identity is the tool that transfers Hessian bounds from a smooth function on a neighborhood of the sphere to its $1$-homogeneous extension. The correction term is the {\em Euler defect} $E$, which vanishes exactly when the function is itself $1$-homogeneous.

\begin{lemma}[Hessian of a homogeneous extension]
\label{Lem:HomExt}
Let $U: \R^d \to \R$ be $C^2$ in a neighborhood of $\mathbb S^{d-1}$ and let $u(\V x)=\|\V x\|_2\,U(\V x/\|\V x\|_2)$ be the $1$-homogeneous extension of its restriction to the sphere. Then, for every $\V \xi\in\mathbb S^{d-1}$ and with $\M P_{\V \xi}=\M I-\V \xi\V \xi^\Tr$,
\begin{align}
\label{Eq:HomExtHessian}
\Op H\{u\}(\V \xi)=\M P_{\V \xi}\,\Op H\{U\}(\V \xi)\,\M P_{\V \xi}+E(\V \xi)\,\M P_{\V \xi}, \qquad E(\V \xi)=U(\V \xi)-\langle\V \xi,\nabla U(\V \xi)\rangle.
\end{align}
In particular, $\|\Op H\{u\}(\V \xi)\|_{S_1}\le\|\Op H\{U\}(\V \xi)\|_{S_1}+(d-1)\,|E(\V \xi)|$, and $E\equiv0$ when $U$ is itself $1$-homogeneous.
\end{lemma}

\begin{proof}
For $\V x\ne\V 0$, let $\V \xi=\V x/\|\V x\|_2\in\mathbb S^{d-1}$. Since $\nabla\V \xi=\M P_{\V \xi}/\|\V x\|_2$, the chain rule gives
\begin{align*}
\nabla u(\V x)=\V \xi\,U(\V \xi)+\M P_{\V \xi}\nabla U(\V \xi)=\V g(\V \xi),
\end{align*}
where $\V g(\V \eta):=\V \eta\,U(\V \eta)+(\M I-\V \eta\V \eta^\Tr)\nabla U(\V \eta)$ is defined by the same formula in a neighborhood of the sphere. A direct differentiation of $\V g$ yields its Jacobian matrix $\Op D\V g=\big[\partial_j g_i\big]_{i,j=1}^d$,
\begin{align*}
\Op D\V g(\V \eta)&=U\M I+\V \eta\nabla U^\Tr+\Op H\{U\}-\langle\V \eta,\nabla U\rangle\M I-\V \eta\nabla U^\Tr-\V \eta\V \eta^\Tr\Op H\{U\}\\
&=E(\V \eta)\,\M I+(\M I-\V \eta\V \eta^\Tr)\,\Op H\{U\}(\V \eta),
\end{align*}
where all the functions are evaluated at $\V \eta$. The chain rule applied once more gives $\Op H\{u\}(\V x)=\Op D\V g(\V \xi)\,\M P_{\V \xi}/\|\V x\|_2$, which reduces to \eqref{Eq:HomExtHessian} when $\|\V x\|_2=1$; that is, when $\V x=\V \xi$. 
The norm bound follows from the two properties of the orthogonal projector $\M P_{\V \xi}$, whose singular values are $1$, with multiplicity $d-1$ (on the tangent space $\V \xi^\perp$), and $0$ (on $\V \xi$): first, $\|\M P_{\V \xi}\M W\M P_{\V \xi}\|_{S_1}\le\|\M W\|_{S_1}$, by the H\"older-type inequality $\|\M A\M W\M B\|_{S_1}\le\|\M A\|_{S_\infty}\|\M W\|_{S_1}\|\M B\|_{S_\infty}$ and $\|\M P_{\V \xi}\|_{S_\infty}=1$; second, $\|\M P_{\V \xi}\|_{S_1}=d-1$, the sum of its singular values. Finally, if $U$ is $1$-homogeneous, then Euler's identity $\langle\V \xi,\nabla U(\V \xi)\rangle=U(\V \xi)$ gives $E\equiv0$.
\end{proof}

\begin{theorem}[$L_1$-bound under a bounded spherical Hessian variation]
\label{Theo:HTVApprox}
Let $q:\R^d\to\R$ be continuous and $1$-homogeneous with $\HTV_{\mathbb S^{d-1}}(q)<\infty$, and let $(\Spc F_n)$ be a sequence of shape-regular and quasi-uniform simplicial fans with maximal spherical cell diameter $h_n\le\frac12$. Then, the preaveraged conic spline \eqref{Eq:SplineAveraged} satisfies
\begin{align}
\label{Eq:ApproxBoundL1}
\|q-\tilde q_n\|_{L_1(\mathbb{S}^{d-1})}\le C''_{\Spc{F}}\,h_n^2\, \HTV_{\mathbb S^{d-1}}(q),
\end{align}
where the constant $C''_{\Spc F}$ depends only on the shape-regularity and quasi-uniformity constants of the fans and on the kernel $\chi$. In particular, $h_n\asymp n^{-1/(d-1)}$ and
\begin{align}
\label{Eq:ApproxRateL1}
\|q-\tilde q_n\|_{L_1(\mathbb{S}^{d-1})} \le C'''_{\Spc{F}}\, \frac{\HTV_{\mathbb S^{d-1}}(q)}{n^{2/(d-1)}}.
\end{align}
\end{theorem}

\begin{proof}
We first recall the reduction to a flat simplex from the proof of Theorem~\ref{Theo:ConicSplineApprox}. Let $C_\alpha=\Cone(\M r_{\alpha,1},\ldots,\M r_{\alpha,d})$ be a cell of $\Spc F_n$ with simplex $T_\alpha$. For any $\V \xi\in C_\alpha\cap\mathbb S^{d-1}$, there exist unique coefficients $a_1,\ldots,a_d\ge 0$ such that $\V \xi=\sum_{j=1}^{d}a_j\M{r}_{\alpha,j}$. Since $\|\V \xi\|_2=1$, not all coefficients vanish, and we set
\begin{align*}
\varsigma=\sum_{j=1}^{d}a_j>0, \qquad \V z=\frac{\V \xi}{\varsigma} =\sum_{j=1}^{d}\frac{a_j}{\varsigma}\,\M{r}_{\alpha,j}\in T_\alpha.
\end{align*}
By construction, a conic spline with nodal values $g(\M r_i)$ is, on $C_\alpha$, the positively homogeneous extension of the affine interpolant $\Op I_{T_\alpha}\{g\}$ of Lemma~\ref{Lem:AffineInterp}. For any $1$-homogeneous $g$, we therefore have that
\begin{align}
\label{Eq:FlatReduction}
g(\V \xi)-\sum_{i=1}^n g(\M r_i)\varphi_i(\V \xi) = \varsigma\bigl(g(\V z)-\Op I_{T_\alpha}\{g\}(\V z)\bigr).
\end{align}
Shape-regularity provides three constants that do not depend on $\alpha$ or $n$: $C_1$ with $\diam(T_\alpha)\le C_1 h_n$; $c>0$ with $\|\V z\|_2\ge c$ for all $\V z\in T_\alpha$, because the vertices of $T_\alpha$ lie on the unit sphere;
and $\varsigma_{\Spc F}$ of Proposition~\ref{Prop:ConicBsplines} with $\varsigma\le\varsigma_{\Spc F}$. We also note that the radial projection $\V z\mapsto\V z/\|\V z\|_2$ maps $T_\alpha$ onto the spherical cell $\Sigma_\alpha=C_\alpha\cap\mathbb S^{d-1}$.

\medskip
\emph{Homogeneous mollification}:
Let $h=h_n\le\frac12$ and let $q_h$ denote the $1$-homogeneous extension of the restriction of $q*\chi_h$ to the unit sphere,
\begin{align}
\label{Eq:HomMollifier}
q_h(\V x)=\|\V x\|_2\,(q*\chi_h)\Big(\frac{\V x}{\|\V x\|_2}\Big)=\int_{\R^d}\chi_h(\V y)\,q\big(\V x-\|\V x\|_2\,\V y\big)\dint\V y,
\end{align}
where the second expression follows from the homogeneity of $q$. The function $q_h$ is $C^\infty$ on $\R^d\backslash\{\V 0\}$, it is $1$-homogeneous, and its nodal values are those of $\tilde q_n$, so that $\tilde q_n$ is the sampled conic spline of $q_h$. We split the error as
\begin{align*}
q-\tilde q_n=(q-q_h)+(q_h-\tilde q_n)
\end{align*}
and bound the two terms separately: the first one is a mollification error and the second one is the interpolation error of a smooth function. Beforehand, we record a consequence of \eqref{Eq:ShellMass}. For $\V \xi\in\mathbb S^{d-1}$ and $0<t\le\frac12$, the ball $B(\V \xi,t)$ is contained in the conic shell $K_{\Sigma(\V \xi,2t)}(1-t,1+t)$, where $\Sigma(\V \xi,\tau)=\{\V \eta\in\mathbb S^{d-1}:\|\V \eta-\V \xi\|_2\le\tau\}$ is the spherical cap of radius $\tau$ around $\V \xi$; hence,
\begin{align}
\label{Eq:BallMass}
|\Op H\{q\}|\big(B(\V \xi,t)\big)\le c_0\,t\,\mu_q\big(\Sigma(\V \xi,2t)\big), \qquad c_0=2\cdot(\tfrac32)^{d-2}.
\end{align}

\medskip
\emph{Mollification error}:
For any $1$-homogeneous function $g$, the change to polar coordinates gives $\|g\|_{L_1(A)}=\varkappa_d\,\|g\|_{L_1(\mathbb S^{d-1})}$ over the shell $A=A(1,2)$, where $A(a,b)=\{\V x: a<\|\V x\|_2<b\}$, with $\varkappa_d=\int_1^2 r^{d}\dint r$. Since $q-q_h$ is $1$-homogeneous, it suffices to bound $\|q-q_h\|_{L_1(A)}$. The substitution $\V y=h\V u$ in \eqref{Eq:HomMollifier} and the evenness of $\chi$ give
\begin{align*}
q(\V x)-q_h(\V x)=\frac12\int_{\R^d}\chi(\V u)\Big(2q(\V x)-q\big(\V x+\V w\big)-q\big(\V x-\V w\big)\Big)\dint\V u,
\end{align*}
where $\V w=h\|\V x\|_2\,\V u$ depends on $\V x$ and $\V u$. We write $\Delta^2_{\V w}q(\V x)=2q(\V x)-q(\V x+\V w)-q(\V x-\V w)$ for the second difference. Let us first assume that $q$ is $C^2$ away from the origin. Taylor's formula with integral remainder, applied to $t\mapsto q(\V x+t\V w)$ on $[-1,1]$, yields
\begin{align*}
\big|\Delta^2_{\V w}q(\V x)\big|\le\|\V w\|_2^2\int_{-1}^1(1-|t|)\,\|\Op H\{q\}(\V x+t\V w)\|_{S_1}\dint t.
\end{align*}
For fixed $\V u$ with $\|\V u\|_2\le1$ and $t\in[-1,1]$, the map $\Phi(\V x)=\V x+th\|\V x\|_2\V u$ has the Jacobian matrix $\M I+th\,\V u\,\V \xi^\Tr$ with $\V \xi=\V x/\|\V x\|_2$, whose determinant $1+th\langle\V u,\V \xi\rangle$ lies in $[1-h,1+h]$. Since $\Phi$ is a perturbation of the identity by a map of Lipschitz constant at most $h<1$, it is injective, and it maps $A$ into $A(1-h,2+2h)$. The integration over $\V x\in A$, where $\|\V w\|_2^2\le 4h^2$, and the change of variables $\V v=\Phi(\V x)$ then give
\begin{align*}
\int_A\big|\Delta^2_{\V w}q(\V x)\big|\dint\V x
&\le \frac{4h^2}{1-h}\,\Big(\int_{-1}^1(1-|t|)\dint t\Big)\,|\Op H\{q\}|\big(A(1-h,2+2h)\big)\nonumber\\
&\le 8h^2\,|\Op H\{q\}|\big(A(\tfrac12,3)\big),
\end{align*}
uniformly in $\V u$. This inequality extends to a general $q$ with bounded spherical Hessian variation by a standard approximation argument: one applies it to the mollified functions $q*\chi_\epsilon$, whose Hessian variation on a set is bounded by that of $q$ on its $\epsilon$-neighborhood. By letting $\epsilon\downarrow0$, the left-hand side converges by the uniform continuity of $q$ on compact sets and the boundary spheres of the shell being $|\Op H\{q\}|$-null by \eqref{Eq:ShellMass}. The integration against $\chi(\V u)$ and \eqref{Eq:ShellMass} finally yield
\begin{align}
\label{Eq:MollificationError}
\|q-q_h\|_{L_1(\mathbb S^{d-1})}=\frac{\|q-q_h\|_{L_1(A)}}{\varkappa_d}&\le \frac{4h^2}{\varkappa_d}\,|\Op H\{q\}|\big(A(\tfrac12,3)\big)\nonumber\\
&\le C_3\,h^2\,\HTV_{\mathbb S^{d-1}}(q).
\end{align}

\medskip
\emph{Hessian of the homogeneous mollification}:
We now bound $\Op H\{q_h\}$ pointwise on the sphere. Lemma~\ref{Lem:HomExt}, applied to $U=q*\chi_h$, gives for $\V \xi\in\mathbb S^{d-1}$
\begin{align*}
\|\Op H\{q_h\}(\V \xi)\|_{S_1}\le\|\Op H\{q*\chi_h\}(\V \xi)\|_{S_1}+(d-1)\,|E_h(\V \xi)|,
\end{align*}
where $E_h(\V \xi)=(q*\chi_h)(\V \xi)-\langle\V \xi,\nabla(q*\chi_h)(\V \xi)\rangle$ is the Euler defect of $q*\chi_h$.
The first term is the convolution of the Hessian measure with the kernel, $\Op H\{q*\chi_h\}=\Op H\{q\}*\chi_h$, so that \eqref{Eq:BallMass} gives
\begin{align*}
\|\Op H\{q*\chi_h\}(\V \xi)\|_{S_1}\le\|\chi_h\|_{L_\infty}\,|\Op H\{q\}|\big(B(\V \xi,h)\big)\le c_0\|\chi\|_{L_\infty}\,h^{1-d}\,\mu_q\big(\Sigma(\V \xi,2h)\big).
\end{align*}
For the second term, we note that $\nabla q$ is a function of locally bounded variation on $\R^d\backslash\{\V 0\}$ that satisfies Euler's identity $q(\V v)=\langle\V v,\nabla q(\V v)\rangle$ for almost every $\V v$. Hence,
\begin{align*}
E_h(\V \xi)&=\int_{\R^d}\chi_h(\V y)\Big(q(\V \xi-\V y)-\langle\V \xi,\nabla q(\V \xi-\V y)\rangle\Big)\dint\V y\\
&=-\int_{\R^d}\chi_h(\V y)\,\langle\V y,\nabla q(\V \xi-\V y)\rangle\dint\V y\\
&=-\frac12\int_{\R^d}\chi_h(\V y)\,\big\langle\V y,\nabla q(\V \xi-\V y)-\nabla q(\V \xi+\V y)\big\rangle\dint\V y,
\end{align*}
where the last step uses the evenness of $\chi$. Let $\bar{\V g}$ denote the mean of $\nabla q$ over $B(\V \xi,h)$. We insert $\pm\bar{\V g}$ and invoke the Poincar\'e inequality $\|\V g-\bar{\V g}\|_{L_1(B)}\le C_P\,h\,|\Op D\V g|(B)$ for functions of bounded variation on a ball $B$ of radius $h$ \cite{Ambrosio2000}, which yields
\begin{align*}
|E_h(\V \xi)|&\le h\,\|\chi_h\|_{L_\infty}\,\|\nabla q-\bar{\V g}\|_{L_1(B(\V \xi,h))}\\
&\le C_P\|\chi\|_{L_\infty}\,h^{2-d}\,|\Op H\{q\}|\big(B(\V \xi,h)\big)\\
&\le c_0C_P\|\chi\|_{L_\infty}\,h^{3-d}\,\mu_q\big(\Sigma(\V \xi,2h)\big).
\end{align*}
Since $h\le\frac12$, the two contributions combine into
\begin{align}
\label{Eq:HessianMollified}
\|\Op H\{q_h\}(\V \xi)\|_{S_1}\le C_4\,h^{1-d}\,\mu_q\big(\Sigma(\V \xi,2h)\big), \qquad \forall\V \xi\in\mathbb S^{d-1}.
\end{align}

\medskip
\emph{Interpolation error of the mollified function}:
We apply \eqref{Eq:FlatReduction} with $g=q_h$, Lemma~\ref{Lem:AffineInterp} with $p=1$, and the homogeneity \eqref{Eq:HessianHomogeneity} of $\Op H\{q_h\}$ together with $\|\V y\|_2\ge c$ on $T_\alpha$, exactly as in the proof of Theorem~\ref{Theo:ConicSplineApprox}. Since the radial projection of $T_\alpha$ is $\Sigma_\alpha$, \eqref{Eq:HessianMollified} yields
\begin{align*}
\|q_h-\tilde q_n\|_{L_\infty(\Sigma_\alpha)}&\le \frac{\varsigma_{\Spc F}C_1^2}{2c}\,h^2\sup_{\V \eta\in\Sigma_\alpha}\|\Op H\{q_h\}(\V \eta)\|_{S_1}\\
&\le \frac{\varsigma_{\Spc F}C_1^2C_4}{2c}\,h^{3-d}\,\mu_q\big(\Sigma_\alpha^{+2h}\big),
\end{align*}
where $\Sigma_\alpha^{+\tau}$ denotes the $\tau$-neighborhood of $\Sigma_\alpha$ in $\mathbb S^{d-1}$. Shape-regularity gives $\Spc H^{d-1}(\Sigma_\alpha)\le c_2h^{d-1}$, so that $\|q_h-\tilde q_n\|_{L_1(\Sigma_\alpha)}\le C_5\,h^2\,\mu_q(\Sigma_\alpha^{+2h})$. Finally, shape-regularity and quasi-uniformity imply that every point of the sphere belongs to at most $N_0$ of the enlarged cells $\Sigma_\alpha^{+2h}$, with $N_0$ independent of $n$. The summation over the cells therefore gives
\begin{align}
\label{Eq:InterpolationErrorMollified}
\|q_h-\tilde q_n\|_{L_1(\mathbb S^{d-1})}\le C_5\,h^2\sum_\alpha\mu_q\big(\Sigma_\alpha^{+2h}\big)\le C_5N_0\,h^2\,\HTV_{\mathbb S^{d-1}}(q).
\end{align}

\medskip
\emph{Conclusion}:
The triangle inequality, \eqref{Eq:MollificationError}, and \eqref{Eq:InterpolationErrorMollified} give \eqref{Eq:ApproxBoundL1} with $C''_{\Spc F}=C_3+C_5N_0$. Since the rays are quasi-uniform on $\mathbb S^{d-1}$, the classical covering estimate for the sphere gives $h_n\asymp n^{-1/(d-1)}$, and the substitution of this relation into \eqref{Eq:ApproxBoundL1} yields \eqref{Eq:ApproxRateL1}.
\end{proof}

Theorems~\ref{Theo:ConicSplineApprox} and~\ref{Theo:HTVApprox} exhibit the expected pairing between the regularity index and the error norm. The index $\HLip_{\mathbb S^{d-1}}$ of \eqref{Eq:HLip} is built on the Schatten-$\infty$ norm and on a supremum over the sphere, and it controls the $L_\infty$ error at the second order. 
Analogously, its companion $\HTV_{\mathbb S^{d-1}}$, which is built on the dual Schatten-$1$ norm and on a measure over the sphere, controls the $L_1$ error at the same second order. The second regime is the one that matters for us: it contains the polyhedral quasi-seminorms themselves, whose Hessian is a measure carried by the facets of the fan, as well as every quasi-seminorm whose spherical profile has bounded Hessian variation, without any differentiability assumption.

Three remarks on the hypotheses are in order. First, no second-order regularity at all is needed for convergence without a rate: if $q$ is merely continuous, then \eqref{Eq:FlatReduction} and the fact that the barycentric coordinates are non-negative and sum to one give $|q(\V z)-\Op I_{T_\alpha}\{q\}(\V z)|\le\omega_q(C_1h_n)$, where $\omega_q$ is the modulus of continuity of $q$ on the compact shell $\{c\le\|\V x\|_2\le1\}$, so that $\|q-q_n\|_{L_\infty(\mathbb S^{d-1})}\le \varsigma_{\Spc F}\,\omega_q(C_1h_n)\to0$ for the sampled spline; the $\ell_p$ quasi-norms of Appendix~\ref{App:ellp} are covered by this statement, and by nothing stronger. Second, when $q$ is itself a polyhedral quasi-seminorm with fan $\Spc G$ and $\Spc F_n$ refines $\Spc G$, the sampled spline is exact, $q_n=q$, since both are linear on each cell of $\Spc F_n$ and agree at its rays. Third, the preaveraging of the nodal values in \eqref{Eq:SplineAveraged} is not a technical convenience but a necessity when $d\ge4$: the spherical profile of $q$ may then exhibit a spike $\V \xi\mapsto\min\{M,\|\V \xi-\V \xi_0\|_2^{-a}\}$ with $0<a<d-3$, whose Hessian variation stays bounded as $M\to\infty$ while the nodal value at a ray $\M r_i=\V \xi_0$ grows like $M$, so that no estimate of the sampled spline in terms of $\HTV_{\mathbb S^{d-1}}(q)$ alone can hold. This is the failure of the embedding of $W^{2,1}$ into $L_\infty$ in dimension $d-1\ge3$. In the low dimensions $d\in\{2,3\}$, where the embedding holds \cite{Demengel1984}, the sampled spline satisfies the same bound; the planar case is made explicit in Proposition~\ref{Prop:PlanarHTVApprox}, with sharp constants. Preaveraging is the price to pay for a statement that is uniform in the dimension.

\section{The Planar Case}
\label{Sec:Planar}

The case $d=2$ deserves a separate treatment. Everything becomes explicit there: the conic B-splines on a uniform polar grid reduce to the shifts of a single trigonometric hat function, and the regularity indices of Section~\ref{Sec:Regularity} reduce to classical quantities of the one-dimensional spline literature.

\subsection{The Angular Profile of a Planar Quasi-seminorm}
\label{Sec:PlanarSetting}
We write a point of $\R^2\backslash\{\V 0\}$ in polar coordinates as $\V x=r\,\V \xi(\theta)$ with $r=\|\V x\|_2>0$, $\V \xi(\theta)=(\cos\theta,\sin\theta)$, and $\theta\in\T=\R/2\pi\Z$, and we denote by $\V \xi'(\theta)=(-\sin\theta,\cos\theta)$ the unit tangent vector to the circle. 

By the homogeneous extension \eqref{Eq:HomogeneousExtension}, a continuous CH functional $q$ on $\R^2$ is entirely encoded by its {\em angular profile}
\begin{align}
\label{Eq:AngularProfile}
\varrho:\T\to[0,\infty), \qquad \varrho(\theta)=q\big(\V \xi(\theta)\big), \qquad q(\V x)=r\,\varrho(\theta).
\end{align}
When $q$ is a quasi-norm, $\varrho=1/r_{B_q}\circ\V \xi$ is the reciprocal of the radial function \eqref{Eq:RadialProfile} of its unit ball: the polar curve $r=\varrho(\theta)$ is the radial graph of $q$, whereas $\partial B_q$ is the polar curve $r=1/\varrho(\theta)$.
Conversely, {\em every} continuous $\varrho:\T\to[0,\infty)$ that is $\pi$-periodic, $\varrho(\theta+\pi)=\varrho(\theta)$, defines through \eqref{Eq:AngularProfile} a functional $q$ of class $\CH(\R^2\to\R_+)$: the continuity and the absolute homogeneity of $q$ are built into the formula, and the $\pi$-periodicity is the transcription of the symmetry $q(-\V x)=q(\V x)$. The unit ball of $q$ is the star-shaped region that is bounded by the polar curve $r=1/\varrho(\theta)$, in accordance with \eqref{Eq:RadialProfile}. Whether $q$ is a quasi-seminorm is decided by the zero set of $\varrho$, and the answer is the planar transcription of Theorem~\ref{Theo:SemiQuasiNorm}, Item 3.

\begin{proposition}[Planar criterion]
\label{Prop:PlanarCriterion}
Let $\varrho:\T\to[0,\infty)$ be continuous and $\pi$-periodic and let $q(\V x)=r\,\varrho(\theta)$ be its homogeneous extension \eqref{Eq:AngularProfile}. Then, $q$ is a quasi-seminorm on $\R^2$ if and only if one of the following three cases occurs.
\begin{enumerate}
\item The profile is strictly positive. Then, $q$ is a quasi-norm with $\Spc N_q=\{\V 0\}$ and with the modulus $C\le\max_\theta\varrho/\min_\theta\varrho$.
\item The profile vanishes at exactly one antipodal pair $\{\theta_0,\theta_0+\pi\}$ and there exist $0<m_0\le M_0<\infty$ such that
\begin{align}
\label{Eq:PlanarFirstOrder}
\forall \theta \in \T: \quad m_0\,|\sin(\theta-\theta_0)| \le \varrho(\theta) \le M_0\,|\sin(\theta-\theta_0)|.
\end{align}
Then, $\Spc N_q=\R\V \xi(\theta_0)$.
\item The profile vanishes identically.
\end{enumerate}
In particular, a profile that vanishes at two or more antipodal pairs without vanishing identically, or that vanishes to an order different from one, defines a CH functional that is {\em not} a quasi-seminorm.
\end{proposition}

\begin{proof}
The null space of $q$ is the union of the lines $\R\V \xi(\theta_0)$ over the zeros $\theta_0$ of $\varrho$, and a union of lines is a linear subspace of $\R^2$ only if it is $\{\V 0\}$, a single line, or $\R^2$. This gives the trichotomy. In the first and third cases, the conclusion follows from Theorem~\ref{Theo:QuasiNormAuto} and from $q\equiv0$, respectively. In the second case, $\dist\big(r\V \xi(\theta),\R\V \xi(\theta_0)\big)=r\,|\sin(\theta-\theta_0)|$, so that \eqref{Eq:PlanarFirstOrder} is exactly the comparison \eqref{Eq:DistanceComparison} of Theorem~\ref{Theo:SemiQuasiNorm}, Item 3.
\end{proof}

The differential quantities of Section~\ref{Sec:RegularityIndices} also have a transparent expression in terms of the profile. If $\varrho$ is twice differentiable, then the differentiation of $q(\V x)=r\,\varrho(\theta)$ in polar coordinates gives
\begin{align}
\label{Eq:PolarGradient}
\nabla q(\V x)&=\varrho(\theta)\,\V \xi(\theta)+\varrho'(\theta)\,\V \xi'(\theta),\\
\label{Eq:PolarHessian}
\Op H\{q\}(\V x)&=\frac{\varrho(\theta)+\varrho''(\theta)}{r}\;\V \xi'(\theta)\,\V \xi'(\theta)^\Tr,
\end{align}
so that the Hessian of a planar $1$-homogeneous function has rank one, with $\V \xi'(\theta)$ as its only eigendirection. The three regularity indices \eqref{Eq:Lipq}, \eqref{Eq:HLip}, and \eqref{Eq:HTVq} therefore read
\begin{align}
\label{Eq:PlanarIndices}
\Lip_{\mathbb S^1}(q)&=\sup_{\theta}\sqrt{\varrho(\theta)^2+\varrho'(\theta)^2}, \qquad
\HLip_{\mathbb S^1}(q)=\sup_{\theta}\,|\varrho(\theta)+\varrho''(\theta)|, \nonumber\\
\HTV_{\mathbb S^1}(q)&=\int_0^{2\pi}|\varrho(\theta)+\varrho''(\theta)|\,\dint\theta,
\end{align}
with the understanding that the last quantity is the total-variation norm $\|\Op L\varrho\|_{\Spc M(\T)}$ of the measure $\Op L\varrho$ when $\varrho$ is merely continuous with $\Op L\varrho=\varrho''+\varrho$ a finite measure. The second-order differential operator
\begin{align}
\label{Eq:OperatorL}
\Op L=\Op D^2+\Op I
\end{align}
is thus the natural regularization operator of the planar problem. Its null space $\spn\{\cos\theta,\sin\theta\}$ consists precisely of the profiles of the linear functions $\V x\mapsto\langle\M f,\V x\rangle$, which are the functions with $\HTV_{\mathbb S^1}=0$; among the $\pi$-periodic profiles, however, $\Op L$ is injective, since $\cos$ and $\sin$ are not $\pi$-periodic. Two further consequences of \eqref{Eq:PolarHessian} are worth recording. First, since the Hessian is rank one, $q$ is convex if and only if $\varrho+\varrho''\ge0$, so that
\begin{align}
\label{Eq:PlanarConvexity}
q \text{ is a seminorm} \iff \Op L\varrho\ge0 \text{ (as a measure)},
\end{align}
in which case $\HTV_{\mathbb S^1}(q)=\int_0^{2\pi}(\varrho+\varrho'')\dint\theta=\int_0^{2\pi}\varrho(\theta)\dint\theta$ is simply the integral of $q$ over the unit circle, because $\int_0^{2\pi}\varrho''(\theta)\dint\theta=\varrho'(2\pi)-\varrho'(0)=0$ by periodicity; when $\Op L\varrho$ is merely a measure, the same conclusion follows from $\langle\Op D^2\varrho,1\rangle=\langle\varrho,\Op D^2 1\rangle=0$. Second, the first-order index $\TV_{\mathbb S^1}(q)=\int_0^{2\pi}\sqrt{\varrho^2+\varrho'^2}\,\dint\theta$
is the length of the polar curve $r=\varrho(\theta)$, the radial graph of $q$---not of $\partial B_q$, whose length is $\int_0^{2\pi}\sqrt{\varrho^2+\varrho'^2}\,\varrho^{-2}\dint\theta$.

The operator \eqref{Eq:OperatorL} is the planar instance of a differential operator that is intrinsic to the sphere. For a $1$-homogeneous $q(r\V \xi)=r\varrho(\V \xi)$ on $\R^d$, the expression of the Laplacian in hyperspherical coordinates, $\Delta=\partial_r^2+\frac{d-1}r\partial_r+\frac1{r^2}\Delta_{\mathbb S^{d-1}}$ \cite[Chapter~IV]{SteinWeiss1971}, \cite[Chapter~3]{AtkinsonHan2012}, gives $\Delta q(r \V \xi)=\frac1r\big(\Delta_{\mathbb S^{d-1}}\varrho(\V \xi)+(d-1)\varrho(\V \xi)\big)$, where $\Delta_{\mathbb S^{d-1}}$ is the Laplace--Beltrami operator of the sphere, so that the natural regularization operator in dimension $d$ is
\begin{align}
\label{Eq:OperatorLd}
\Op L_d=\Delta_{\mathbb S^{d-1}}+(d-1)\Op I,
\end{align}
of which $\Op L=\Op D^2+\Op I$ is the case $d=2$.
Its null space is the eigenspace of $\Delta_{\mathbb S^{d-1}}$ that is associated with the eigenvalue $-(d-1)$: the spherical harmonics of degree $\ell$ are the eigenfunctions of $\Delta_{\mathbb S^{d-1}}$ with the eigenvalue $-\ell(\ell+d-2)$ \cite[Chapter~3]{Groemer1996}, \cite[Chapter~3]{AtkinsonHan2012}, and the degree $\ell=1$ consists of the restrictions to the sphere of the linear functions, in exact analogy with $\spn\{\cos,\sin\}$. Since $\Delta_{\mathbb S^{d-1}}\varrho+(d-1)\varrho$ is the trace of the Hessian of $q$ on the sphere, the identity $\HTV_{\mathbb S^{d-1}}(q)=(d-1)\int_{\mathbb S^{d-1}}\varrho\,\dint\Spc H^{d-1}$ holds for every convex member of the class in every dimension, which generalizes the identity that follows \eqref{Eq:PlanarConvexity}. What is specific to the plane is that the Hessian \eqref{Eq:PolarHessian} has rank one, so that its nuclear norm coincides with its trace and the convexity criterion is a scalar inequality.
\subsection{Polyhedral Quasi-seminorms on a Uniform Polar Grid}
\label{Sec:PlanarPolyhedral}
A complete simplicial fan of $\R^2$ is specified by its rays $\V \xi(\theta_1),\dots,\V \xi(\theta_{N'})$ with $\theta_1<\dots<\theta_{N'}<\theta_1+2\pi$, the cells being the sectors between consecutive rays; the simpliciality of the cells requires that each sector be narrower than $\pi$. The fan is symmetric when its rays come in antipodal pairs, in which case $N'=2N$ is even. We now focus on the {\em uniform polar grid} with $2N$ rays,
\begin{align}
\label{Eq:UniformGrid}
\theta_n=nh, \qquad h=\frac{\pi}{N}, \qquad n=0,\dots,2N-1, \qquad N\ge2.
\end{align}
On the sector $\theta\in[nh,(n+1)h]$, a conic spline is linear, $q(\V x)=\langle\M f_n,\V x\rangle$, so that its profile is $\varrho(\theta)=\langle\M f_n,\V \xi(\theta)\rangle=a_n\cos\theta+b_n\sin\theta$. In other words, {\em the profile of a planar polyhedral quasi-seminorm is a continuous, piecewise-trigonometric function of the first order}; that is, a trigonometric spline of order $2$ in the sense of Schoenberg \cite{Schoenberg1964,Lyche1979}, whose knots are the rays of the fan. The conic B-spline \eqref{Eq:ConicBspline} that is associated with the ray $\theta_n$ is then $\varphi_n(\V x)=r\,\phi(\theta-nh)$, where $\phi$ is the $2\pi$-periodic extension of the {\em trigonometric hat function}
\begin{align}
\label{Eq:TrigHat}
\phi(\theta)=
\begin{cases}
\dfrac{\sin(h-|\theta|)}{\sin h}, & |\theta|\le h,\\[1.5ex]
0, & h<|\theta|\le\pi.
\end{cases}
\end{align}
Indeed, for $\theta\in[nh,(n+1)h]$, the function $\theta\mapsto\sin\big((n+1)h-\theta\big)/\sin h$ is of the form $a\cos\theta+b\sin\theta$, it takes the value $1$ at $\theta=nh$ and the value $0$ at $\theta=(n+1)h$, which is the biorthogonality condition \eqref{Eq:Biorthogonality}; the same holds on the sector to the left. The uniform grid makes all the conic B-splines shifts of a single generator, exactly as for the polynomial B-splines on a uniform grid. The function $\phi$ is non-negative, it vanishes at the other nodes, and the expansion \eqref{Eq:qBsplineExpansion} of a polyhedral CH functional becomes
\begin{align}
\label{Eq:PlanarExpansion}
\varrho(\theta)=\sum_{n=0}^{2N-1}c_n\,\phi(\theta-nh), \qquad c_n=\varrho(nh)\ge0, \qquad c_{n+N}=c_n,
\end{align}
so that $q$ is parametrized by the $N$ free non-negative numbers $\V c=(c_0,\dots,c_{N-1})$. The reproduction of the linear functions on any fan translates into the identities $\sum_n\cos(nh)\,\phi(\theta-nh)=\cos\theta$ and $\sum_n\sin(nh)\,\phi(\theta-nh)=\sin\theta$. By Proposition~\ref{Prop:PlanarCriterion}, the parametrized functional is a quasi-norm whenever $\V c>\V 0$, and
it is a quasi-seminorm if and only if at most one of the free parameters $c_0,\dots,c_{N-1}$ vanishes; the first-order condition \eqref{Eq:PlanarFirstOrder} is automatic here, because $\varrho$ is piecewise trigonometric. When $q$ is a quasi-norm, the boundary vertices \eqref{Eq:Boundaryvertices} of the star-shaped polygon $B_q$ are the points $\V \xi(nh)/c_n$.

The regularity indices are now evaluated from the nodal values. On the sector $[nh,(n+1)h]$, \eqref{Eq:PlanarExpansion} gives $\varrho(\theta)=\big(c_n\sin((n+1)h-\theta)+c_{n+1}\sin(\theta-nh)\big)/\sin h$, whose derivative at the two endpoints is
\begin{align*}
\varrho'(nh^+)=\frac{c_{n+1}-c_n\cos h}{\sin h}, \qquad \varrho'\big((n+1)h^-\big)=\frac{c_{n+1}\cos h-c_n}{\sin h}.
\end{align*}
Since $\Op L\varrho=\varrho''+\varrho$ vanishes inside each sector, the measure $\Op L\varrho$ is a sum of Dirac impulses at the knots, whose weights are the jumps of $\varrho'$:
\begin{align}
\label{Eq:PlanarLrho}
\Op L\varrho=\sum_{n=0}^{2N-1}\frac{c_{n-1}-2\cos h\;c_n+c_{n+1}}{\sin h}\;\delta(\cdot-nh)=\sum_{n=0}^{2N-1}\frac{[\M L\V c]_n}{\sin h}\;\delta(\cdot-nh),
\end{align}
where the indices are taken modulo $2N$ and where $\M L$ is the circulant matrix of the {\em trigonometric second difference} with stencil $(1,\,-2\cos h,\,1)$. This matrix annihilates the sequences $\cos(nh)$ and $\sin(nh)$, as it should, and it reduces to the ordinary second difference $\Delta^2$ up to a term of order $h^2$. By comparing \eqref{Eq:PlanarLrho} with \eqref{Eq:PlanarIndices}, and by invoking the symmetry $c_{n+N}=c_n$, we obtain the planar version of \eqref{Eq:HTVCPWL},
\begin{align}
\label{Eq:PlanarHTV}
\HTV_{\mathbb S^1}(q)=\frac{1}{\sin h}\sum_{n=0}^{2N-1}\big|c_{n-1}-2\cos h\;c_n+c_{n+1}\big|=\frac{2}{\sin h}\,\|\M L\V c\|_{\ell_1},
\end{align}
where, in the last expression, $\V c\in\R^N$ collects the free parameters and $\M L\in\R^{N\times N}$ is the corresponding circulant matrix with indices modulo $N$. Likewise, the convexity criterion \eqref{Eq:PlanarConvexity} becomes the componentwise condition $\M L\V c\ge\V 0$,
which is the planar transcription of the vertex condition of a polyhedral norm: the node $\V \xi(nh)/c_n$ is an extreme point of $B_q$ when $[\M L\V c]_n>0$, and it lies in the interior of an edge when $[\M L\V c]_n=0$.

The reduced matrix $\M L\in\R^{N\times N}$ is diagonalized by the discrete Fourier transform, with the eigenvalues
\begin{align}
\label{Eq:PlanarEigen}
\lambda_k=2\Big(\cos\big(\tfrac{2\pi k}{N}\big)-\cos\big(\tfrac{\pi}{N}\big)\Big), \qquad k=0,\dots,N-1,
\end{align}
none of which vanishes, since $\cos(2\pi k/N)=\cos(\pi/N)$ would require $2k=\pm1$ modulo $2N$. The operator $\M L$ is therefore invertible on the $\pi$-periodic sequences, which is the discrete counterpart of the injectivity of $\Op L$ that was observed in Section~\ref{Sec:PlanarSetting}: the sequences $\cos(nh)$ and $\sin(nh)$ that it annihilates are $2\pi$-periodic and change sign under $n\mapsto n+N$, so that they are not of the form \eqref{Eq:PlanarExpansion}. The practical consequences are that the penalty \eqref{Eq:PlanarHTV} is a norm of $\V c$---whence the coercivity of a learning problem that uses it as a penalty---and that the condition number of $\M L$, which grows like $N^2$, is the quantity that governs the convergence of the splitting algorithm.

The planar case is also where the preaveraging of the nodal values in Theorem~\ref{Theo:HTVApprox} can be dispensed with, and where the constants become explicit.

\begin{proposition}[Approximation in the plane under a bounded Hessian variation]
\label{Prop:PlanarHTVApprox}
Let $q$ be a continuous quasi-seminorm on $\R^2$ with profile $\varrho$ and $\HTV_{\mathbb S^1}(q)=\|\Op L\varrho\|_{\Spc M(\T)}<\infty$, and let $q_N$ be its sampled conic spline \eqref{Eq:SplineInterpolant} on the uniform polar grid \eqref{Eq:UniformGrid} with $2N$ rays and $h=\pi/N\le\pi/2$. Then, with $I_n=\big(nh,(n+1)h\big)$,
\begin{align}
\label{Eq:PlanarHTVApproxLinf}
\|q-q_N\|_{L_\infty(\mathbb S^1)}&\le \tfrac12\tan\big(\tfrac h2\big)\,\max_{n}\,|\Op L\varrho|(I_n)\ \le\ \frac{h}{\pi}\,\HTV_{\mathbb S^1}(q),\\
\label{Eq:PlanarHTVApproxL1}
\|q-q_N\|_{L_1(\mathbb S^1)}&\le \tfrac12\,h\tan\big(\tfrac h2\big)\,\HTV_{\mathbb S^1}(q)\ \le\ \frac{h^2}{\pi}\,\HTV_{\mathbb S^1}(q).
\end{align}
In particular, if $q$ is polyhedral with knots on the grid, then $q_N=q$.
\end{proposition}

\begin{proof}
Let $\varrho_N$ be the profile of $q_N$. On the sector $I_n$, $\varrho_N$ is of the form $a\cos\theta+b\sin\theta$, since $q_N$ is linear in Cartesian coordinates there, and it agrees with $\varrho$ at the two endpoints; because $h<\pi$, it is the unique solution of $\Op L\varrho_N=0$ on $I_n$ with these boundary values. The error $e=\varrho-\varrho_N$ therefore solves $\Op L e=\Op L\varrho$ on $I_n$ with $e(nh)=e\big((n+1)h\big)=0$. The Green's function of $\Op L=\Op D^2+\Op I$ on an interval of length $h<\pi$ with homogeneous Dirichlet conditions is
\begin{align*}
G(\theta,\tau)=-\frac{\sin\big(\min(\theta,\tau)-nh\big)\,\sin\big((n+1)h-\max(\theta,\tau)\big)}{\sin h}
\end{align*}
for $\theta,\tau\in[nh,(n+1)h]$, which is continuous, vanishes when $\tau$ is an endpoint, and satisfies $|G|\le\sin^2(h/2)/\sin h=\frac12\tan(h/2)$, the maximum being attained at $\theta=\tau=(n+\frac12)h$. Since $\Op L\varrho$ is a finite measure and $G$ is continuous, we have the representation $e(\theta)=\int_{I_n}G(\theta,\tau)\,\dint(\Op L\varrho)(\tau)$, in which the atoms of $\Op L\varrho$ at the endpoints contribute nothing. Hence, $\sup_{I_n}|e|\le\frac12\tan(h/2)\,|\Op L\varrho|(I_n)$. The maximum over $n$ gives the first bound of \eqref{Eq:PlanarHTVApproxLinf}, and the integration over $I_n$ followed by the summation over $n$ gives the first bound of \eqref{Eq:PlanarHTVApproxL1}, because the $2N$ open sectors are disjoint. The second bounds follow from $\tan x\le\frac4\pi x$ on $[0,\frac\pi4]$, which is the convexity of the tangent there, applied with $x=h/2$. Finally, if $q$ is polyhedral with knots on the grid, then $\Op L\varrho$ is carried by the grid points, so that $e\equiv0$.
\end{proof}

The two bounds are the planar instances of Theorems~\ref{Theo:ConicSplineApprox} and~\ref{Theo:HTVApprox}, with the sampled spline in both and with the pairing between index and error norm made explicit: the mass of $\Op L\varrho$ on a single sector controls the uniform error, its total mass controls the mean error, and the two coincide, up to the factor $h$, when that mass is spread evenly.

\section{Conclusion}
\label{Sec:Conclusion}
The three structural requirements that we have imposed on a regularizer---continuity, symmetry, and amplitude equivariance---select the functionals $g=q^\gamma$ with $q$ continuous and absolutely homogeneous, and convexity is the only thing that separates the seminorms of \cite{UnserDucotterd2026} from the much larger CH class. Within that class, the conic splines play the role that the polyhedral norms play in the convex theory:
they are dense, by Theorems~\ref{Theo:ConicSplineApprox} and~\ref{Theo:HTVApprox}, and they are parametrized by their nodal values along the rays of a fan, with no constraint other than the non-negativity and the symmetry of those values.

Two features of the non-convex setting deserve to be remembered. The first is favorable: a positive-definite CH functional is a quasi-norm for free, with a modulus that is the ratio of the largest to the smallest value on the unit sphere, so that the admissibility of a learned functional needs no monitoring as long as its nodal values stay positive. The second is a genuine obstruction: a null space cannot be obtained by letting the nodal values degenerate, because the quasi-triangle inequality then requires the null set to be a linear subspace along which the functional vanishes at the first order. The construction to be used instead is the analysis form of Theorem~\ref{Theo:QuasiSeminormStructure}, Item 1(iii), where a linear map carries the invariances and the conic spline carries the positive-definite factor.

The natural continuations of this work are the learning of a planar quasi-seminorm from samples, which turns into a sparse regression problem in the nodal values with the weighted $\ell_1$-norm \eqref{Eq:PlanarHTV} as the penalty, its extension to $d\ge3$ through the operator $\Op L_d=\Delta_{\mathbb S^{d-1}}+(d-1)\Op I$ of \eqref{Eq:OperatorLd} and a symmetric triangulation of the sphere,
and the integration of the resulting parametrization into amplitude-equivariant neural architectures, where the conic spline would play the role of the final pooling layer.
\appendix
\section*{APPENDICES}
\section{Quasi-norms in Infinite Dimensions}
\label{App:InfiniteDim}
The next theorem provides a purely algebraic characterization of quasi-norms as gauge functionals on an arbitrary vector space $\Spc X$. Unlike Theorem~\ref{Theo:GaugeSemiQuasiNorm}, no topology is involved, so that the result applies equally well in infinite dimensions. As a consequence, no regularity can be inferred: the quasi-norm $q$ is not necessarily continuous, and the associated set $B$ need not be closed, in contrast with the finer characterization of Theorem~\ref{Theo:GaugeSemiQuasiNorm} for $\Spc X=\R^d$. The result is classical; see \cite[pp. 159--160]{Koethe1983}.

\begin{theorem}[Quasi-norms as quasi-gauges]
\label{Theo:QuasiNormsGauges}
Let $\Spc X$ be a vector space. A functional $q:\Spc X\to[0,\infty)$ is a quasi-norm if and only if there exists a set $B\subseteq \Spc X$ such that
\begin{enumerate}
\item $B$ is balanced;
\item $B$ is absorbing;
\item $B$ is star-shaped with respect to the origin;
\item there exists $a_0\geq 2$ such that $B+B\subseteq a_0 B$;
\item $\bigcap_{\lambda>0}\lambda B=\{0\}$;
\end{enumerate}
and
\begin{align*}
q(x)=\Gamma_B(x):=\inf\{\lambda>0:\ x\in \lambda B\}.
\end{align*}
Moreover, if $q$ is given, then one may take $B=B_q:=\{ x\in\Spc X:\ q(x)\leq 1\}$.
\end{theorem}

\begin{proof}
\emph{Necessity}:
Let us assume that $q$ is a quasi-norm and let $B = B_q = \{x \in \Spc X : q(x) \leq 1\}$.
\begin{itemize}
    \item \textit{Balanced}: if $x \in B$ and $|\alpha| \leq 1$, then $q(\alpha x) = |\alpha| \, q(x) \leq 1$, so that $\alpha x \in B$.
    \item \textit{Absorbing}: let $x \in \Spc X$. If $x = 0$, then $x \in \lambda B$ for every $\lambda>0$. Otherwise, positive-definiteness gives $q(x)>0$ and $q\big(x/q(x)\big) = 1$, so that $x \in q(x)B$.
    \item \textit{Star-shaped}: this follows directly from $B$ being balanced.
    \item \textit{Inclusion property}: by the quasi-triangle inequality, $q(x+y) \leq C(q(x)+q(y)) \le 2C$ for $x,y \in B$. Thus, $x+y \in 2C B$ and the property holds with $a_0 = 2C$.
    \item \textit{Definiteness}: $\bigcap_{\lambda>0}\lambda B = \{x : q(x) \leq \lambda \ \forall \lambda > 0\} = \{x : q(x) = 0\} = \{0\}$.
\end{itemize}
Finally, for every $x\in\Spc X$ and $\lambda>0$, we have that $x\in \lambda B_q \iff q(x/\lambda)\leq 1 \iff q(x)\leq \lambda$, so that
\begin{align*}
\Gamma_{B_q}(x) =\inf\{\lambda>0:\ q(x)\leq \lambda\} =q(x).
\end{align*}

\medskip
\emph{Sufficiency}:
Let us assume that $B$ satisfies (1)--(5) and let $q = \Gamma_B$. Since $B$ is absorbing, $q(x)<\infty$ for every $x$.
\begin{itemize}
    \item \textit{Homogeneity}: since $B$ is balanced, $x \in \lambda B \iff \alpha x \in |\alpha|\lambda B$ for $\alpha \ne 0$, so that $q(\alpha x) = |\alpha| \, q(x)$. In particular, $q(0)=0$.
    \item \textit{Definiteness}: if $q(x)=0$, then, for every $\lambda>0$, there exists $\mu<\lambda$ such that $x\in \mu B$. Since $B$ is star-shaped, $\mu B\subseteq \lambda B$ and hence $x\in\lambda B$ for every $\lambda>0$. Property (5) then gives $x=0$.
    \item \textit{Quasi-triangle inequality}: let $x, y \in \Spc X$. For any $\epsilon > 0$, we have that $x \in \big(q(x)+\epsilon\big)B$ and $y \in \big(q(y)+\epsilon\big)B$. Let $\rho_\epsilon = \max\{q(x)+\epsilon, q(y)+\epsilon\}$. Since $B$ is star-shaped, $x, y \in \rho_\epsilon B$, so that $x+y \in \rho_\epsilon (B+B) \subseteq \rho_\epsilon a_0 B$. Hence, $q(x+y) \leq a_0 \rho_\epsilon$ and, by letting $\epsilon \downarrow 0$, $q(x+y) \le a_0 \max\{q(x), q(y)\} \leq a_0\big(q(x)+q(y)\big)$.
\end{itemize}
Thus, $q$ is a quasi-norm.
\end{proof}

When $\Spc X=\R^d$, Conditions (2) and (5) may be replaced by
\begin{align*}
\V 0\in \intr(B) \qquad\text{and}\qquad B \text{ bounded}.
\end{align*}
Indeed, in finite dimensions, $\V 0\in\intr(B)$ is equivalent to $B$ being absorbing. Moreover, if $B$ is bounded, then $\bigcap_{\lambda>0}\lambda B=\{\V 0\}$ since, for every $\V x\neq \V 0$, the choice $\lambda<\|\V x\|/R$ with $B\subseteq R\,B_{\|\cdot\|}$ gives $\V x\notin \lambda B$.

\section{Counterexamples}
\label{App:Counterexamples}
The examples of this appendix mark the boundary between the convex theory and its non-convex extension. Each of them exhibits a property that every (semi-)norm on $\R^d$ enjoys and that a member of the CH class may lack, and each has a consequence for the way the theory is set up.

\subsection{A Quasi-norm That Is Not Continuous}
\label{App:Counterexample}
By contrast with the triangle inequality, the quasi-triangle inequality does not imply continuity.

\begin{example}[A discontinuous quasi-norm on $\R^2$]
\label{Ex:Discontinuous}
Let $\Spc L=\R\V e_1$ and let $q:\R^2\to[0,\infty)$ be given by
\begin{align}
\label{Eq:Counterexample}
q(\V x)=
\begin{cases}
\|\V x\|_2, & \V x\in \Spc L,\\[1ex]
2\|\V x\|_2, & \V x\in \R^2\backslash \Spc L;
\end{cases}
\end{align}
that is, $q$ is the homogeneous extension \eqref{Eq:HomogeneousExtension} of the profile that equals $1$ at $\pm\V e_1$ and $2$ elsewhere on $\mathbb S^1$. Since $\Spc L$ is a linear subspace, $\V x$ and $\lambda\V x$ fall under the same branch of \eqref{Eq:Counterexample}, so that $q$ is absolutely homogeneous; it is positive-definite; and $\|\V x\|_2\le q(\V x)\le2\|\V x\|_2$ gives $q(\V x+\V y)\le2\|\V x+\V y\|_2\le2\big(q(\V x)+q(\V y)\big)$. Hence, $q$ is a quasi-norm with $C=2$. Yet it is discontinuous at every point of $\Spc L\backslash\{\V 0\}$: for $\V y_n=(1,\frac1n)\to\V e_1$, we have that $q(\V y_n)\to2\ne1=q(\V e_1)$. Its unit ball is the disk of radius $\frac12$ with two antennas of length $1$ along $\Spc L$; this set is symmetric, closed, bounded, star-shaped, with the origin in its interior, and it satisfies $B_q+B_q\subseteq4B_q$, so that the only hypothesis of Theorem~\ref{Theo:GaugeSemiQuasiNorm}, Item 2,
that fails is the continuity of the radial function.
\end{example}

\emph{Implications}: The algebraic conditions of Theorem~\ref{Theo:QuasiNormsGauges} do not deliver continuity, which therefore has to be assumed; this is why it is built into the CH class (Definition~\ref{Def:CHscalar}) and into the gauge theorem. The requirement is not a topological one---the functional above is still equivalent to the Euclidean norm---but an approximation-theoretic one: the conic splines of Section~\ref{Sec:PolyhedralQuasiNorms} are continuous, and a uniform limit of continuous functions is continuous, so that a discontinuous quasi-norm is beyond the reach of Theorem~\ref{Theo:ConicSplineApprox} and of any parametrization by continuous functions.

\subsection[Non-Lipschitz Quasi-norms: The \texorpdfstring{$\ell_p$}{lp} Quasi-norm]{Quasi-norms That Are Not Lipschitz-continuous: The \texorpdfstring{$\ell_p$}{lp} Quasi-norm}
\label{App:ellp}
All norms on $\R^d$ are Lipschitz-continuous, by the equivalence of the norms. The prototypical quasi-norm is continuous but not Lipschitz-continuous.

\begin{proposition}[The $\ell_p$ quasi-norm, $0<p<1$]
\label{Prop:ellp}
Let $\|\V x\|_p:=\big(\sum_{n=1}^d |x_n|^p\big)^{1/p}$ on $\R^d$, with $0<p<1$. Then, the following holds.
\begin{enumerate}
\item $\|\cdot\|_p$ satisfies the $p$-triangle inequality $\|\V x+\V y\|_p^p \le \|\V x\|_p^p+\|\V y\|_p^p$.
\item $\|\cdot\|_p$ is a continuous quasi-norm with modulus $C_p = 2^{\frac1p-1}$, which is sharp and independent of $d$; it is not a norm when $d\ge 2$.
\item When $d\ge2$, $\|\cdot\|_p$ is not Lipschitz-continuous on any neighborhood of a point $\V x\ne\V 0$ with a vanishing coordinate.
\end{enumerate}
\end{proposition}

\begin{proof}
Item 1 follows from $|x_n+y_n|^p\le(|x_n|+|y_n|)^p\le|x_n|^p+|y_n|^p$, which is the subadditivity of $t\mapsto t^p$ on $\R_+$, summed over $n$. For Item 2, the absolute homogeneity and the positive-definiteness are read off the definition; the $p$-triangle inequality gives $\big|\,\|\V x+\V h\|_p^p-\|\V x\|_p^p\big|\le\|\V h\|_p^p$, hence the continuity, and, once raised to the power $1/p$ and combined with the convexity bound $(a+b)^{1/p}\le2^{\frac1p-1}(a^{1/p}+b^{1/p})$, the quasi-triangle inequality with $C_p=2^{\frac1p-1}$. For $\V x,\V y$ with disjoint supports and $A=\|\V x\|_p^p$, $B=\|\V y\|_p^p$, the quotient $\|\V x+\V y\|_p/(\|\V x\|_p+\|\V y\|_p)=(A+B)^{1/p}/(A^{1/p}+B^{1/p})$ equals $2^{\frac1p-1}$ for $A=B$, so that $C_p$ is sharp; the choice $\V x=\V e_1$, $\V y=\V e_2$ realizes it and shows that the triangle inequality fails, since $2^{1/p}>2$. For Item 3, let $x_n=0$ and $a=\|\V x\|_p^p>0$, so that $\|\V x+t\V e_n\|_p^p=a+|t|^p$. The convexity of $s\mapsto(a+s)^{1/p}$ gives $\|\V x+t\V e_n\|_p-\|\V x\|_p\ge\frac1p a^{\frac1p-1}|t|^p$, so that the difference quotient along $\V e_n$ grows like $|t|^{p-1}$ as $t\to0$.
\end{proof}

Item 3 has a geometric reading: the gauge of a star body is Lipschitz-continuous if and only if the body is star-shaped with respect to a ball centered at the origin \cite{Boulkhemair2023}, and the unit ball $B_p$ is star-shaped with respect to the origin only. Indeed, for $\V y=\varepsilon\V e_2$ and $s=1-\delta$,
\begin{align*}
\|(1-s)\V y+s\V e_1\|_p^p=s^p+(1-s)^p\varepsilon^p=1+\delta^p\varepsilon^p-p\,\delta+O(\delta^2)>1
\end{align*}
for $\delta$ small, because $\delta^p\gg\delta$; the segment $[\V y,\V e_1]$ leaves $B_p$ near $\V e_1$, so that no ball $B(\V 0,r)$ is contained in $\ker(B_p)$, however small $\varepsilon$.
The boundary of $B_p$ reaches $\pm\V e_n$ tangentially to the axes, at cusps, so that no cone with apex at the origin and axis $\pm\V e_n$ fits inside the ball near those points.

\emph{Implications}: The Lipschitz index \eqref{Eq:Lipq} of Section~\ref{Sec:Regularity} is infinite for $\|\cdot\|_p$, so that the approximation theorems of Section~\ref{Sec:Regularity} give no rate. The conic spline interpolant still converges uniformly on the sphere, by the uniform continuity of the profile, and the exact rate is dictated by the cusps along the coordinate axes.

\begin{proposition}[Saturation at low smoothness]
\label{Prop:ellpRate}
Let $d=2$ and $0<p<1$, let $q=\|\cdot\|_p$, and let $q_N$ be the conic spline interpolant of $q$ on the uniform polar grid \eqref{Eq:UniformGrid} with $2N$ rays and $h=\pi/N$. Then, there exist two constants $0<c\le C<\infty$ that depend on $p$ alone and are such that
\begin{align}
\label{Eq:ellpRate}
c\,h^{p}\ \le\ \|q-q_N\|_{L_\infty(\mathbb S^{1})}\ \le\ C\,h^{p}.
\end{align}
\end{proposition}

\begin{proof}
The angular profile is $\varrho(\theta)=\big(|\cos\theta|^p+|\sin\theta|^p\big)^{1/p}$, which satisfies $\varrho(\theta)=1+\frac1p|\theta|^p+O(\theta^2)$ as $\theta\to0$ and is $C^2$ away from the multiples of $\pi/2$. On a cell that contains no multiple of $\pi/2$, the interpolation error is $O(h^2\sup|\varrho''|)$, the supremum being taken over that cell; the summation of these contributions over the cells at distance $kh$ from an axis gives a total that is $O(h^2)+O(h^p)$, the second term coming from the cells adjacent to an axis. On a cell adjacent to an axis, $[0,h]$ say, the chord of $\varrho$ differs from $\varrho$ by $\frac1p\big(2^{-p}-\tfrac12\big)h^p+O(h^2)$ at the midpoint $h/2$, which is bounded below by $c\,h^p$ for $h$ small because $2^{-p}>\tfrac12$ for $p<1$. The two estimates give \eqref{Eq:ellpRate}.
\end{proof}

The exponent $p<1$ in \eqref{Eq:ellpRate} is the exact counterpart of the exponent $2$ in \eqref{Eq:ApproxRate}: the conic spline interpolant converges at the order dictated by the smoothness of the profile, and the second-order rate of Theorem~\ref{Theo:ConicSplineApprox} is an upper limit that only the twice-differentiable members of the class attain.
The cusps of $B_p$ are also what makes $\|\cdot\|_p$ a sparsity-promoting penalty, so that the loss of regularity is the price of the sparsity. Finally, the free bound $C\le M/m$ of Theorem~\ref{Theo:QuasiNormAuto}, applied with the Euclidean norm, returns $d^{\frac1p-\frac12}$ here, against the sharp $C_p=2^{\frac1p-1}$: the bound guarantees the quasi-triangle inequality, but it can be very pessimistic about its modulus.

\subsection{CH Functionals That Are Not Quasi-seminorms}
\label{App:NotQuasiSeminorm}
A positive-definite CH functional is a quasi-norm by Theorem~\ref{Theo:QuasiNormAuto}. When the null space is non-trivial, the quasi-triangle inequality holds if and only if the null space is a linear subspace and $q$ is comparable to the distance to it (Theorem~\ref{Theo:SemiQuasiNorm}, Item 3); each of the two requirements can fail on its own.

\begin{example}[A CH functional whose null space is not a linear subspace]
\label{Ex:NotQuasiSeminorm}
The conic spline $q(\V x)=\big|\,|x_1|-|x_2|\,\big|$ on $\R^2$ is linear on each of the eight sectors delimited by the axes and the diagonals, and its null space is the union of the two diagonals. It is not a quasi-seminorm, since $q(1,1)=q(1,-1)=0$ while $q\big((1,1)+(1,-1)\big)=q(2,0)=2$.
\end{example}

\begin{example}[A CH functional that is not comparable to the distance to its null space]
\label{Ex:NotQuasiSeminorm2}
The functional $q(\V x)=\sqrt{\|\V x\|_2\,|x_2|}$ on $\R^2$ has the linear null space $\R\V e_1$, but $q(t,1)=(1+t^2)^{1/4}$ is unbounded although $\dist\big((t,1),\R\V e_1\big)=1$: the upper bound of \eqref{Eq:DistanceComparison} fails. So does the quasi-triangle inequality, since $q\big((t,0)+(0,1)\big)=q(t,1)$ is unbounded while $q(t,0)+q(0,1)=1$.
\end{example}

\emph{Implications}: The quasi-triangle inequality has to be checked whenever a null space is present, and the check has two parts. The first one is combinatorial and concerns the zero set: for a conic spline, an isolated vanishing nodal value, as in Example~\ref{Ex:NotQuasiSeminorm}, is what Proposition~\ref{Prop:ConicSplineAdmissible} forbids. The second one is metric and cannot be read off the null space alone: Example~\ref{Ex:NotQuasiSeminorm2} has the right null space and still fails, because $q$ vanishes too fast---like $|x_2|^{1/2}$---in the direction that leaves it; in the plane, it is the first-order vanishing of the profile that Proposition~\ref{Prop:PlanarCriterion} requires. The safe construction is the invariant case of Theorem~\ref{Theo:QuasiSeminormStructure}, where a linear map carries the null space and the non-convex factor stays positive-definite. None of this arises in the convex setting: if $q\in\CH(\R^d\to\R_+)$ is convex, then $q(\V x+\V y)=2\,q\big(\tfrac12\V x+\tfrac12\V y\big)\le q(\V x)+q(\V y)$, so that a convex CH functional is a seminorm, and both examples above are non-convex---the first one because a convex non-negative functional that vanishes on two lines vanishes on their convex hull, the second one because it is not even constant along its null space ($q(0,1)=1$, $q(3,1)=10^{1/4}$).
\section{Complements on Polyhedral Geometry}
\label{App:Polyhedral}
\begin{proposition}[Simplicial refinement without new rays]
\label{Prop:FanSimplicial}
Let $\Spc F$ be a finite polyhedral fan in $\R^d$ whose cones are pointed. Then, there exists a simplicial fan $\widetilde{\Spc F}$ such that $|\widetilde{\Spc F}|=|\Spc F|$ and such that, for every $\widetilde C\in\widetilde{\Spc F}$, there exists $C\in\Spc F$ with $\widetilde C\subseteq C$. Moreover, $\widetilde{\Spc F}$ may then be chosen so that its rays are among the rays of $\Spc F$. For a general fan, the refinement still exists, but new rays may be needed.
\end{proposition}

\begin{proof}
Let $\Spc R(\Spc F)$ denote the finite set of rays of $\Spc F$. We intersect the fan with the unit sphere:
\begin{align*}
\Spc S:=\{C\cap \mathbb S^{d-1}:\ C\in\Spc F\}.
\end{align*}
This is a finite spherical polyhedral complex whose vertices are induced by the rays of $\Spc F$.
By the standard triangulation theorem for finite polyhedral complexes, each cell of $\Spc S$ admits a triangulation that uses only its existing vertices, and these triangulations can be chosen face-to-face; see, e.g., \cite[Chap.~5]{Ziegler1995}. The coning of each simplex of this triangulation from the origin gives a collection of simplicial cones. The resulting fan $\widetilde{\Spc F}$ satisfies $|\widetilde{\Spc F}|=|\Spc F|$, and every one of its cones is contained in a cone of $\Spc F$. Moreover, no new ray is introduced, since the vertices of the spherical triangulation are precisely induced by the rays of $\Spc F$.
\end{proof}

Let $q:\R^d\to\R$ be linear on each cone of a complete polyhedral fan $\Spc F$. If $\widetilde{\Spc F}$ is a simplicial refinement of $\Spc F$, then $q$ is linear on each cone of $\widetilde{\Spc F}$. Indeed, if $q(\V x)=\langle \M f_C,\V x\rangle$ for $\V x\in C$, then $q(\V x)=\langle \M f_C,\V x\rangle$ for every $\widetilde C\subseteq C$ and every $\V x\in \widetilde C$.

\begin{definition}[Continuous piecewise-linear function]
\label{Def:CPWL}
A function $f: \R^d \to \R$ is \emph{continuous piecewise-linear} (CPWL) if
\begin{enumerate}
\item it is globally continuous, with $f \in C(\R^d)$, and
\item there exists a finite collection of closed domains $\{C_\alpha\}_{\alpha\in A}$ whose union is $\R^d$ and such that the restriction of $f$ to each $C_\alpha$ is affine:
\begin{align*}
f|_{C_\alpha}(\V x) = \M a_\alpha^\Tr \V x + b_\alpha,
\end{align*}
where $\M a_\alpha \in \R^d$ and $b_\alpha \in \R$.
\end{enumerate}
\end{definition}

\begin{proposition}[Polyhedral partition and face restriction]
\label{Prop:PartitionCPWL}
Let $f: \R^d \to \R$ be a CPWL function.
\begin{enumerate}
    \item There exists a finite polyhedral complex $\Spc{P}$ that partitions $\R^d$ and whose cells $C_\alpha$ are convex polyhedra.
    \item For each cell $C_\alpha \in \Spc{P}$, the restriction $f|_{C_\alpha}$ is affine, with $f|_{C_\alpha}(\V x) = \M a_\alpha^\Tr \V x + b_\alpha$ for some $\M a_\alpha \in \R^d$ and $b_\alpha \in \R$.
    \item If $F$ is a $k$-dimensional face of $C_\alpha$ with $0 \le k < d$, then the restriction of $f$ to $F$ is an affine function from $\R^k$ to $\R$.
\end{enumerate}
\end{proposition}

\begin{proof}
Let $\{g_{\alpha}\}_{\alpha \in I}$ be the finite set of distinct affine functions that represent the local pieces of $f$, with $g_\alpha(\V x)=\M a_\alpha^\Tr \V x + b_\alpha$. Let $\Spc A$ be the arrangement of the finitely many hyperplanes $H_{\alpha,\beta}=\{\V x\in\R^d: g_\alpha(\V x)=g_\beta(\V x)\}$ over the pairs $\alpha\ne\beta$ for which $g_\alpha-g_\beta$ is not constant. Its closed cells are convex polyhedra, and they form a polyhedral complex that covers $\R^d$. On the interior $\Omega$ of a cell, the functions $g_\alpha$ take pairwise distinct values at every point, so that the sets $\{\V x\in\Omega: f(\V x)=g_\alpha(\V x)\}$ are closed in $\Omega$, pairwise disjoint, and cover $\Omega$. Since $\Omega$ is connected, only one of them is non-empty; hence, $f$ is affine on $\Omega$ and, by continuity, on its closure. We take the cells of $\Spc A$ as the cells $C_\alpha$ of the complex $\Spc P$; they may be finer than the maximal regions of linearity of $f$, which is harmless. (The closures of the connected components of the set on which $f$ is differentiable are not convex in general: for $f(\V x)=\max\{0,\min\{x_1,x_2\}\}$, the component on which $f=0$ is the non-convex set $\{x_1<0\}\cup\{x_2<0\}$; see \cite{Scholtes2012} for the structure theory of piecewise-affine functions.)

Let now $C_\alpha \in \Spc{P}$ be a convex polyhedron that is defined by the intersection of $m$ closed half-spaces, and let $F$ be a $k$-dimensional face of $C_\alpha$. By definition, $F$ is obtained by turning a valid subset of the bounding inequalities into equalities, so that $F = C_\alpha \cap H$ where $H$ is an affine subspace of dimension $k$. Any point $\V x \in F$ is parametrized by an anchor point $\V x_0 \in F$ and a basis matrix $\M{U} \in \R^{d \times k}$ whose columns span the linear subspace parallel to $H$, with $\V x = \V x_0 + \M{U}\V y$ for $\V y \in \R^k$. The evaluation of $f$ on $F$ then yields
\begin{align*}
f|_F(\V y) &= \M a_\alpha^\Tr (\V x_0 + \M{U}\V y) + b_\alpha \\
&= (\M a_\alpha^\Tr \M{U})\V y + (\M a_\alpha^\Tr \V x_0 + b_\alpha),
\end{align*}
which is affine in $\V y$.
\end{proof}

\subsection*{Acknowledgments}
The author used Claude (Anthropic) as an interactive assistant during the preparation of
this paper, for proofreading, for checking and simplifying proofs. All results were verified by the author,
who takes full responsibility for the content. This research was supported by the Swiss National Science Foundation (SNSF) under Grant 200020\_219356.

\bibliographystyle{siam}
\bibliography{QuasiNorm}

\end{document}